\documentclass{article}
\usepackage[british]{babel}
\usepackage{color, amsmath, amssymb, amsthm, mathrsfs, bm, graphicx, color, hyperref, calc}
\usepackage[a4paper,body={16cm,23cm}]{geometry}
\usepackage[all]{xy}
\usepackage{todonotes}

\definecolor{orange}{rgb}{1,0.64,0}
\definecolor{dblue}{rgb}{0,0,0.7}
\newtheoremstyle{mythm}{1ex}{1ex}{\it\color{dblue}}{}{\bf\color{dblue}}{.}{ }{}
\theoremstyle{mythm}
\newtheorem{thm}{Theorem}
\newtheorem{prop}[thm]{Proposition}

\newtheorem{lem}[thm]{Lemma}

\newlength{\algowidth}
\newenvironment{algo}[1]{%
\vskip\parsep\begin{center}\begin{minipage}{\algowidth}%
\noindent \textbf{#1}%
\begin{list}{}{}%
}{%
\end{list}\end{minipage}\end{center}%
}

\theoremstyle{definition}
\newtheorem*{assumption}{Assumption}
\newtheorem*{defn}{Definition}
\newtheorem*{conj}{Conjecture}

\renewcommand{\geq}{\geqslant}
\renewcommand{\leq}{\leqslant}
\newcommand{\sss}{\scriptscriptstyle}

\DeclareMathOperator{\re}{Re}
\DeclareMathOperator{\im}{Im}

\DeclareMathOperator{\SL}{SL}
\DeclareMathOperator{\GL}{GL}

\DeclareMathOperator{\lcm}{lcm}
\DeclareMathOperator{\tors}{tors}
\DeclareMathOperator{\co}{co}
\DeclareMathOperator{\sign}{sign}

\newcommand{\FF}{\mathbb{F}}
\newcommand{\QQ}{\mathbb{Q}}
\newcommand{\ZZ}{\mathbb{Z}}
\newcommand{\RR}{\mathbb{R}}
\newcommand{\CC}{\mathbb{C}}
\newcommand{\PP}{\mathbb{P}}

\newcommand{\cyclic}[1]{{}^{\ZZ}\!/\!{}_{#1 \ZZ}}
\newcommand{\ve}{\varepsilon}
\newcommand{\no}[1]{\Vert #1\Vert}
\newcommand{\sm}[4]{\bigl(\begin{smallmatrix} #1 & #2 \\ #3 & #4 \end{smallmatrix}\bigr)}
\newcommand{\bigO}{\mathbf{O}}

\DeclareFontFamily{U}{wncy}{}
\DeclareFontShape{U}{wncy}{m}{n}{<->wncyr10}{}
\DeclareSymbolFont{mcy}{U}{wncy}{m}{n}
\DeclareMathSymbol{\Sha}{\mathord}{mcy}{"58}

\begin{document}
\title{Numerical modular symbols for elliptic curves}
\author{Christian Wuthrich}
\maketitle
\abstract{We present a detailed analysis of how to implement the computation of modular symbols for a given elliptic curve by using numerical approximations. This method turns out to be more efficient than current implementations as the conductor of the curve increases.}

\section{Introduction}

  The aim of the article is to describe an alternative algorithm for computing modular symbols for a given fixed elliptic curve~$E/\QQ$ of conductor~$N$. The current implementations use linear algebra with rational coefficients to determine the space of modular symbols attached to~$E$ within the space of all symbols of level~$N$.
  Instead we wish to compute efficiently the value of a modular symbols for a fixed $E$ avoiding to work with the full space whose dimension grows linearly in~$N$. We build on the work of Goldfeld~\cite{goldfeld} using numerical approximation to path integrals in the upper half plane. He already noted that the Atkin-Lehner involutions can be used to avoid integrating close to the real line where the convergence is bad. Goldfeld obtained that a single modular symbol for a semistable elliptic curve can be computed roughly in $N^{1/2}$ steps. See Theorem~\ref{goldfeld_thm} where we recall the precise statement.

  We improve on his work in several directions. First we prove all the rigorous bounds and we present the finer details of an implementation that returns provably correct rational numbers. This uses some theoretical knowledge about the possible denominators. Moreover, we explain what methods can be used for elliptic curves that are not semistable. Furthermore, we explain an idea that allows us to compute the modular symbols at all rational numbers with a fixed denominator. This is very useful for the practical applications we have in mind, for instance computing the $p$-adic $L$-functions of $E$. Finally, we analyse where possible the complexity of our algorithms. For instance, we prove that the set of all Manin symbols can be evaluated in roughly $N^{7/4}$ steps. See Theorem~\ref{complexity_thm} for the precise statement.

  We implemented the algorithms in \texttt{SageMath}~\cite{sage}. It turns out to be faster when computing a single modular symbol and allows for computations with larger conductors than all previous implementations.

  In order to describe the methods and results in more detail, we start by defining modular symbols (sometimes called modular elements). Let $E/\QQ$ be an elliptic curve of conductor $N$. Let $f = \sum a_n q^n$ be the newform of weight $2$ and level $\Gamma_0(N)$ associated to the isogeny class of $E$. We know that $f$ exists by modularity~\cite{modularity}. Given a rational number $r\in\QQ$, we consider the integrals
  \begin{equation}\label{la_eq}
    \lambda(r) = 2\pi i\int_{i\infty}^{r} f(z) dz = 2\pi  \int_{0}^{\infty} f(r+yi) dy \ \in \ \CC.
  \end{equation}
  Let $\gamma^{+}$ be a $\ZZ$-basis of the subgroup of $H_1\bigl(E(\CC),\ZZ\bigr)$ fixed by complex conjugation and, similarly, let $\gamma^{-}$ be a generator for the subgroup on which complex conjugation acts by multiplication with $-1$. Let $\omega_E$ be a N\'eron differential on $E$. Let $\Omega^{+}$ be the smallest positive period of $E$, i.e., $\Omega^{+} = \bigl\vert\int_{\gamma^{+}} \omega_E\bigr\vert$. Similarly, we set $\Omega^{-}$ to be $\bigl\vert \int_{\gamma^{-}}\omega_E\bigr\vert$ in $\RR_{>0}$. The period lattice $\Lambda_E$ of $E$ is either $\ZZ\Omega^{+} \oplus \ZZ\Omega^{-}\,i$ or $\ZZ\Omega^{+}\oplus \ZZ\tfrac{1}{2} (\Omega^{+}+\Omega^{-}\,i)$ depending on whether the discriminant of $E$ is negative or positive.

  Manin~\cite{manin72} and Drinfeld~\cite{drinfeld} showed that there exists an integer $t$ such that $t\cdot \lambda(r) \in \ZZ\Omega^{+} \oplus\ZZ\Omega^{-}\,i$ for all $r\in\QQ$. In Section~\ref{denominator_sec}, we will look for a good bound on $t$ in practice. Now we define the rational numbers
  \begin{equation*}
    [r]^{+} = \frac{\re\bigl(\lambda(r)\bigr)}{\Omega^{+}} \, \in \QQ\qquad\text{ and }\qquad
    [r]^{-} = \frac{\im\bigl(\lambda(r)\bigr)}{\Omega^{-}}\,\in\QQ.
  \end{equation*}
  In this article, the map $r\mapsto [r]^{\pm}$ will be called a modular symbol rather than the homological version where the paths are called modular symbols. Our main goal is to find a fast algorithm for computing the values $[r]^{\pm}$ for a given curve $E$ and $r\in \QQ$.

  Current implementations of modular symbols of $E$ compute the values $[r]^{\pm}$ as follows. First they determine the vector space over $\QQ$ of all modular symbols as the modular form varies through all rational cuspidal forms of weight $2$ and level $N$. Then the matrices for the first few Hecke operators are computed and they are used, together with the known eigenvalues $a_p$ for our curve $E$, to find the subspace corresponding to our fixed cuspform $f$. (Or rather quotient as they work with the dual space.) Once this initial step of finding a basis for this subspace is done, the value of $[r]^{\pm}$ for a given $r$ is computed efficiently using the continued fractions expansion of the rational number $r$.

  A thorough explanation of this method is given in Stein's book~\cite{stein} and in Cremona's book~\cite{cremona}. It is implemented in Cremona's library \texttt{eclib}~\cite{eclib}, \texttt{Magma}~\cite{magmamanual}, \texttt{PARI/GP}~\cite{pari} and~\cite{pariscript}, and \texttt{SageMath}~\cite{sage}. Originally these implementations were written to find the elliptic curves of a given conductor as explained in~\cite{cremona}. In particular, the modularity of $E$ was proven with this method, too.

  Instead, we use here that the modularity of the elliptic curve is known. We wish to avoid to work with the space of all modular symbols of level $N$ because this involves manipulations with sparse matrices of size $N/3\times N/4$ as explained in \S~8.9 of~\cite{stein}. As $N$ increases the initial step takes up a very long time and it currently makes it difficult to work with elliptic curves of conductor larger than $10^5$.

  The approach in this paper is to compute the values of $\lambda(r)\in\CC$ by finding a numerical approximation to the integral in~\eqref{la_eq}. We assume that we are given the values of the Fourier coefficients $a_n$ of $f$; for instance \texttt{PARI}~\cite{pari} yields these very fast by point counting on the reductions of $E$. We also know how to compute good approximations to the values of the periods $\Omega^{\pm}$. We make one assumption: We suppose that the Manin constant of the strong Weil curve in the isogeny class of $E$ is $1$. See Section~\ref{outstanding_subsec} for the concrete implication of this assumption.

  Here is how the modular symbol $[r]^{+}$ is computed in practice. First we use Manin's trick~\cite{manin72} with continued fractions to split the path from $i\infty$ to $r$ into pieces (Section~\ref{manins_trick_sec}). This reduces our problem to evaluating so called Manin symbols (Section~\ref{manin_symbols_subsec}). These are integrals between two cusps $r$ and $r'$. The main advantage is that the denominators of $r$ and $r'$ are now small compared to $N$. The path from $r$ to $r'$ is split up at the best place into two pieces (Section~\ref{unitary_sec}). We use an Atkin-Lehner involution as in~\cite{goldfeld} to move the path close to $r$ to a path close to $i\infty$ where the Fourier expansion of $f$ allows for fast integration. This integration is done by a summation where the number of terms and the precision of the floating point numbers is determined rigorously to guarantee the result within a given error (Section~\ref{num_sec}).

  However, this is not possible for all cusps $r$. A cusp is called ``unitary'' if it is in the orbit of $i\infty$ under the group of Atkin-Lehner involutions. If we encounter a non-unitary cusp, we have to fall back to a much slower method using so-called transportable paths (Section~\ref{non_unitary_sec}), which we would like to avoid, if at all possible. The most important idea for this is to replace the curve by its quadratic twist of minimal conductor (Section~\ref{twists_subsec}).  Furthermore there is also some flexibility in the continued fraction method.

  The main application we have in mind is to compute algebraic $L$-values $L(E,\chi,1)$ for Dirichlet characters $\chi$ and to compute $p$-adic $L$-functions. In both cases one only needs to find all values $[\tfrac{a}{m}]$ for a fixed $m$. Typically they are all unitary symbols. In Section~\ref{partial_sec}, we explain an idea using partial sums that allows us to evaluate all of these symbols almost as fast as a single evaluation. This has also theoretical implications for the complexity estimates proven in Theorem~\ref{complexity_thm}.

  The structure of the paper goes through the above explanation of the computation in reversed order. It is important first to understand the bounds for the possible denominators of $[r]^{\pm}$ in Section~\ref{denominator_sec}. Then we deal with the numerical approximation in Section~\ref{num_sec} followed by how to split up and move the integration paths in Section~\ref{unitary_sec} and~\ref{non_unitary_sec}. How to use and compute Manin symbols is explained in Section~\ref{manins_trick_sec}. Then, Section~\ref{tweaks_sec} describes how to take advantage of quadratic twists and partial sums and Section~\ref{complexity_sec} looks at the complexity of all steps for unitary symbols.

  We end the paper with examples and numerical comparisons with current implementations. We will illustrate that our method proves to be much faster when we need to evaluate a single, or a small number of values of $[r]^{\pm}$.  It is even comparable when the task is to evaluate all Manin symbols as long as we assume that the curve is semistable. When $N$ is really large, say $10^{10}$, our method still determines single values of modular symbols quite fast, while the current implementations cannot perform the initial step any more. We refer to Section~\ref{example_sec} for precise timings.

  The methods in this paper could be extended to modular forms that do not come from elliptic curves; for instance forms associated to $\QQ$-curves. We have not explored this or any potential generalisations to other groups or situations.

\subsection*{Acknowledgements}
  It is a pleasure to thank John Cremona, Christophe Delaunay, Marc Masdeu, Dave Parkin and Fredrik Str\"omberg for help with the research and the implementation.


\section{Denominator of modular symbols}\label{denominator_sec}

  We will compute a numerical approximation to the rational numbers $[r]^{\pm}$ defined in the introduction. In order to know to what precision we need to compute the approximation, we have to find a good bound on the denominator of the rational numbers $[r]^{+}$ and $[r]^{-}$. This will also lead us to the issue concerning the Manin constant.  See~\cite{wuthrich_int_kato} for further investigations on these denominators.

  First, we need a few further definitions. Throughout this text $E$ will be an elliptic curve defined over $\QQ$ of conductor $N$. We know that $E$ is modular and so let
  \begin{equation*}
    \xymatrix@1{\varphi \colon X_0(N) \ar[r] & E }
  \end{equation*}
  be a modular parametrisation of minimal degree sending $i\infty$ to $O$. This is defined up to an automorphism of $E$ defined over $\QQ$, so up to multiplication by $[-1]$. The Manin constant is defined to be the rational number $c_{\sss E}>0$ such that
  \begin{equation*}
    \varphi^{*}(\omega_{\sss E}) = \pm c_{\sss E} \cdot2\pi i \, f(z) dz.
  \end{equation*}
  We choose $\varphi$ uniquely such that it is the $+$ sign that appears in the above equality.

  In the isogeny class of $E$ there is a unique strong Weil curve $E_0$ (also called the $X_0$-optimal curve in~\cite{stevens}).
  \begin{assumption}
    The Manin constant $c_0=c_{E_0}$ of the strong Weil curve $E_0$ is $1$.
  \end{assumption}
  It is known that $c_0$ is an integer~\cite{edixhovn}, and it is believed to be equal to $1$ in all cases. See~\cite{manin_constant} for a discussion of known results about $c_0$. In particular, it is known that $c_0$ is either $1$ or $2$ when $E$ is semistable. See Section~\ref{outstanding_subsec} for an explanation of how harmful the above assumption is. A consequence of this assumption is that the N\'eron lattice $\Lambda_0$ of $E_0$ is equal to the lattice generated by all values $\lambda(r)-\lambda(s)$ as $r$ and $s$ run through all pairs of $\Gamma_0(N)$-equivalent cusps.

  When comparing the period lattices of $E$ and $E_0$, the quotient of the N\'eron periods will become important. Define the rational numbers
  \begin{equation*}
    q^{+} = \frac{\Omega_E^{+}}{\Omega_{E_0}^{+}} \qquad\text{ and }\qquad q^{-} = \frac{\Omega_E^{-}}{\Omega_{E_0}^{-}}.
  \end{equation*}

  Let $r$ be a rational number.
  \begin{defn}
    We write $r=\tfrac{a}{m}$ as a reduced fraction of integers. Let $M$ be the greatest common divisor of $m$ and the conductor $N$. Hence we can write $N = Q\cdot M$ and $m = d \cdot M$. Following Section~3.1 in~\cite{mazur_sd}, we call a cusp $r$ \textit{unitary} if $Q$ and $M$ are coprime. The integer $Q/\gcd(M,Q)$ is called the \textit{width} of the cusp $r$; for unitary cusps it is simply $Q$.
  \end{defn}
  It is known that $r$ is unitary if and only if the cusps $r$ and $i\infty$ on $X_0(N)$ are in the same orbit under the action by the group of Atkin-Lehner involutions. In the application we have in mind, no prime of additive reduction for $E$ divides the denominator $m$. Then the cusp $\tfrac{a}{m}$ is unitary. For a semistable curve all cusps are unitary.

  Further, we set $\delta^2$ to be the largest square dividing $N$. Thus $\delta=1$ if and only if $E$ is semistable.

  \begin{prop}
   Let $E/\QQ$ be an elliptic curve of conductor $N$. Choose a few primes $\ell>2$ coprime to $N$, with $\ell\equiv 1 \pmod{\delta}$ and set $t_0$ to be the greatest common divisor of the number $N_{\ell}$ of points on the reduction of $E$ modulo $\ell$. Let $t^{\pm}$ be the numerator of $t_0 q^{\pm}$. Assume $c_0=1$. We have
   \begin{equation*}
     [r]^{\pm} \in \frac{c_{\infty}(E_0)}{2 t^{\pm}} \, \ZZ,
   \end{equation*}
   where $c_{\infty}(E_0)$ is the number of connected components of $E_0(\RR)$. If $r$ is unitary, we also get
   \begin{equation*}
     [r]^{+} \in \frac{c_{\infty}(E)}{2 \cdot \# E(\QQ)_{\tors}} \qquad\text{ and }\qquad [r]^{-} \in \frac{1}{2}\ZZ.
   \end{equation*}
   where $c_{\infty}(E)$ is the number of connected components of $E(\RR)$.
  \end{prop}
  For a semistable curve, even without assuming $c_0=1$, we get the bounds $4 \cdot\# E(\QQ)_{\tors}/c_{\infty}(E)\leq 24$ and $4$ for the denominator of $[r]^{+}$ and $[r]^{-}$ respectively. If $t$ is a bound for the denominator of a modular symbols $[r]^{+}$ as above then, in our implementation we now round $t \cdot \re\bigl(\lambda(r)\bigr) /\Omega^{+}$ to the closest integer and find $[r]^+$ by dividing again by $t$. Hence we must compute $\re\bigl(\lambda(r)\bigr)$ with a proven error smaller than $\Omega^+/(2t)$.

  \begin{proof}
   Consider the modular parametrisation $\varphi_0\colon X_0(N)\to E_0$. After identifying $E_0(\CC)$ with $\CC/\Lambda_{E_0}$ via the integration of $\omega_{E_0}$, we get an induced map $\tilde\varphi_0$ from the upper half plane to $\CC/\Lambda_{E_0}$. We find
   \begin{equation*}
    \tilde\varphi_0(r)
      \equiv \int\limits_{O}^{\varphi_0(\Gamma_0(N)r)} \omega_{E_0}
      \equiv \int\limits_{\Gamma_0(N)i\infty}^{\Gamma_0(N)r} \varphi_0^*(\omega_{E_0})
      \equiv c_0 \int\limits_{i\infty}^{r} 2\pi i f(z) dz
      \equiv c_0 \lambda(r) \pmod{ \Lambda_{E_0} }
   \end{equation*}
   By the theorem of Manin and Drinfeld, the modular parametrisation $\varphi\colon X_0(N)\to E$ maps the cusp $r$ to the torsion point $\varphi(r)\in E(\bar{\QQ})$. The action of the Galois group on the cusps on $X_0(N)$ is given in Theorem~1.3.1 in~\cite{stevens_book}. The cusps on $X_0(N)$, and hence all points $\varphi(r)$ for $r\in\QQ$, are defined over the cyclotomic field $K = \QQ(\zeta_\delta)$. The image of the unitary cusps is in the torsion subgroup of $E(\QQ)$ instead.

   If $\ell\equiv 1 \pmod{\delta}$, then there is a place $v$ in $K$ above $\ell$ with residue field $\FF_{\ell}$. If $\ell\nmid N$, then we get a reduction map $E_0(K)\to \tilde E_0(\FF_{\ell})$ of elliptic curves. Since $v\mid \ell$ is unramified and $\ell>2$, we conclude from Theorem~VII.3.4 in~\cite{sil1} that the reduction map is injective on torsion points in $E(K)$. Hence $N_{\ell}$ is a multiple of the order of the torsion subgroup of $E(K)$ for all these $\ell$. We conclude that $t_0 \,\varphi_0(r) = O$ in $E_0(K)$. Therefore $t_0\,\tilde\varphi_0(r)$ and hence $t_0\, c_0\, \lambda(r)$ belong to $\Lambda_{E_0}$. Recall that $\lambda(r) = [r]^{+} q^{+} \,\Omega_{E_0}^{+} + [r]^{-} q^{-} \,\Omega_{E_0}^{-} \, i$. Now if $\Lambda_{E_0}$ is rectangular, then $c_{\infty}(E_0)=2$ and $c_0 t_0 q^{\pm} \,[r]^{\pm}\in \ZZ$. If $\Lambda_{E_0}$ is not rectangular, then $c_{\infty}(E_0)=1$ and  $c_0 t_0 q^{\pm} \,[r]^{\pm}\in \tfrac{1}{2}\ZZ$. Thus combined we find that $c_0 t_0 q^{\pm} \,[r]^{\pm}$ belong to $c_{\infty}(E_0)/2\, \ZZ$.

   Finally if $r$ is unitary, then $\varphi(r)$ belongs to $E(\QQ)$ and hence to $E(\RR)$. This implies in both the rectangular and the non-rectangular case that $[r]^{-} \in \tfrac{1}{2c_0}\ZZ$.  We find that $2 \cdot\# E(\QQ)_{\tors} \,[r]^{+}/ (c_0 \,c_{\infty}(E))$ belongs to $\ZZ$.
  \end{proof}

  By the way, the original proof of Manin~\cite{manin72} and Drinfeld~\cite{drinfeld} used the Hecke operators and found that $N_{\ell}$ for $\ell \equiv 1 \pmod{N}$ is a bound for the order of $\varphi_0(r)\in E(K)$; our bound involving $\delta$ rather than $N$ is better.

  We add the example of the strong Weil curve 121d1. Here $E(\QQ)$ is trivial and $E(\RR)$ is connected. So we expect a denominator $1$ or $2$ for all unitary cusps. We have $\delta=11$ for this curve and $N_{23} = 25$ which is also the greatest common divisor of the first few $\ell\equiv 1 \pmod{11}$. In fact the torsion subgroup of $E(\QQ(\mu_{11}))$ is isomorphic to $\cyclic{25}$. Hence we can bound the denominator of $[r]^{\pm}$ by $50$. One can show that $[r]^{+} \in \tfrac{1}{2}\ZZ$ and $[r]^{-}\in \frac{1}{50}\ZZ$. For instance, $\lambda\bigl(\tfrac{3}{11}\bigr) = -\tfrac{1}{2}\Omega^{+} + \tfrac{27}{50}\Omega^{-}\,i$.

\subsection{Implementation of the Manin constant}

 We add here an explanation of how to implement the Manin constant under the above assumption (as it is now done in SageMath).
 \begin{prop}\label{manin_prop}
  Let $E$ be an elliptic curve defined over $\QQ$.
  Let $n^{\pm}$ be the numerator of $q^{\pm}$ as defined above.
  Then the Manin constant $c_E$ is equal
  \begin{equation*}
     c_E  = \begin{cases}
                    c_0\cdot \tfrac{1}{2} \cdot n^{+}\cdot n^{-} &\text{\parbox{77mm}{if $n^{+}$ and $n^{-}$ are both even and \\ $E_0(\RR)$ has more components than $E(\RR)$,}}\\[3ex]
                    c_0\cdot 2\cdot n^{+}\cdot n^{-} &\text{\parbox{77mm}{if $E(\RR)$ has more components than $E_0(\RR)$ and \\at least one of the denominators of $q^{+}$ and $q^{-}$ is odd,}}\\[2ex]
                    c_0\cdot n^{+}\cdot n^{-} &\text{otherwise. }
                  \end{cases}
  \end{equation*}
\end{prop}
If $\alpha\colon E\to E'$ is an isogeny defined over $\QQ$, then set $n_{\alpha}^{\pm}$ to be the numerator and $d_{\alpha}^{\pm}$ the denominator of $\Omega_{E'}^{\pm}/\Omega_{E}^{\pm}$. Further we define $c_{\alpha}$ by $\alpha^{*}(\omega_{E'}) = c_{\alpha} \,\omega_{E}$. We may choose the differentials such that $c_{\alpha}$ is a positive integer.

Let $\psi\colon E_0\to E$ be the isogeny of smallest degree. By the definition of the strong Weil curve, the modular parametrisation of $E$ factors through $\varphi_0$ and $\psi$. Thus the Manin constant $c_E$ is equal to $c_0\cdot c_{\psi}$. It is convenient to prove a lemma first.

\begin{lem}
 Let $\alpha\colon E\to E'$ be a cyclic isogeny defined over $\QQ$ of degree $p^k$ for some prime $p$. If $p=2$ assume that $E(\RR)$ and $E'(\RR)$ have the same number of connected components. Then $\gcd(n_{\alpha}^{+},n_{\alpha}^{-})=\gcd(d_{\alpha}^{+},d_{\alpha}^{-})=1$ and $c_{\alpha}=n_{\alpha}^{+}\cdot n_{\alpha}^{-}$.
\end{lem}
\begin{proof}
  Integrating against the fixed N\'eron differentials $\omega$ and $\omega'$ on $E$ and $E'$ respectively, we identify $E(\CC)$ with $\CC/\Lambda_E$ and $E'(\CC)$ with $\CC/\Lambda_{E'}$. The isogeny $\CC/\Lambda_E \to \CC/\Lambda_{E'}$ is then induced by the multiplication by $c_{\alpha}$ on $\CC$. Recall also that $c_{\alpha}$ divides $p^k$ as $c_{\alpha}\cdot c_{\hat\alpha} = c_{[p^k]}=p^k$.

  If $n_{\alpha}^{+}$ and $n_{\alpha}^{-}$ have a common divisor $n$, then the isogeny $\alpha$ would factor through $[n]\colon\CC/n\Lambda_{E'}\to \CC/\Lambda_{E}$, but that is not possible as $\alpha$ is cyclic. Similarly $d_{\alpha}^{+}$ and $d_{\alpha}^{-}$ are coprime.

  Since $\ker(\alpha)$ is not a direct sum, it is either contained in $E(\RR)$ or in $E(\CC)^{-}$, the set of points $Q$ in $E(\CC)$ whose complex conjugate is $\bar Q =-Q$. Let $z+\Lambda_E$ be a generator of $\ker(\alpha)$. Now if $z=x+iy$, then $\bar z \equiv \pm z \pmod{\Lambda_E}$ implies that either $2x\in\Lambda_E$ or $2iy\in\Lambda_{E}$.

  Assume now that $p$ is odd. Then the above implies that either $\ker(\alpha)$ is generated by $\Omega_{E}^{+}/p^k +\Lambda_E$ or it is generated by $i\Omega_{E}^{-}/p^k +\Lambda_E$. In the first case, we have $c_\alpha \Omega_E^{+}/p^k = \Omega_{E'}^{+}$ and $c_{\alpha} \Omega_E^{-} =\Omega_{E'}^{-}$. Together with $c_{\alpha}\mid p^k$, this implies that $n_{\alpha}^{+}=d_{\alpha}^{-}=1$ and $c_{\alpha}=n_{\alpha}^{-}$. The lemma is then proved in this case. The second case, when $\Omega_E^{-} i/p^k$ is in $\ker(\alpha)$, is similar but with signs swapped.

  Finally, we assume $p=2$. Consider $w = 2^{k-1}z$. Then $w+\Lambda$ is a $2$-torsion point on $E$ and, since $\alpha$ is defined over $\QQ$, it lies in $E(\RR)[2]$. First, if $E(\RR)$ is connected, then $w\in\Omega^{+}_E/2 +\Lambda_E$ and $z \in \Omega^{+}_E/2^k + \Lambda_E$. We find ourselves in a case in which the explanation for general $p$ treated above extends to $p=2$. Also when $E(\RR)$ has two connected components, we fall back onto the two cases treated above, except when $w\in (\Omega_E^{+}+i\Omega_E^{-})/2 +\Lambda_E$. However in this last case, $E'(\RR)$ is connected, which is excluded by assumption.
\end{proof}
\begin{proof}[Proof of Proposition~\ref{manin_prop}]
 We factor $\psi = \beta\circ \alpha$ with $\alpha\colon E_0\to E'$ and $\beta\colon E'\to E$. We can impose that $\alpha$ is cyclic and $E_0(\RR)$ and $E'(\RR)$ have the same number of connected components and that $\beta$ has degree $2$ when $E_0(\RR)$ and $E(\RR)$ do not have the same number of connected components otherwise $\beta$ is trivial.

 Decomposing $\alpha$ into isogenies of prime power degrees, we can apply the previous lemma repeatedly. It follows that $\gcd(n_{\alpha}^{+},n_{\alpha}^{-})=\gcd(d_{\alpha}^{+},d_{\alpha}^{-})=1$ and $c_{\alpha}=n_{\alpha}^{+}\cdot n_{\alpha}^{-}$. This concludes the case when $E_0(\RR)$ and $E(\RR)$ have the same number of connected components.

 Assume now that $E_0(\RR)$ has two and $E(\RR)$ has one connected component. As seen above in the case $p=2$, it follows that the kernel of $\beta$ is generated by $(\Omega_{E'}^{+}+i\Omega_{E'}^{-})/2+\Lambda_{E'}$. Either $c_{\beta}=1$ or $2$. In the first case, we have $\Omega_{E'}^{\pm} = \Omega_{E}^{\pm}$ and hence $c_{\psi}= c_{\alpha}\cdot c_{\beta} = n_{\alpha}^{+}\cdot n_{\alpha}^{-}= n_{\psi}^{+}\cdot n_{\psi}^{-}$ proves this case. In the second case, we have $c_{\beta}=2=\Omega_{E}^{\pm}/\Omega_{E'}^{\pm}$. Thus $n_{\psi}^{\pm}/d_{\psi}^{\pm} = 2 n_{\alpha}^{\pm}/d_{\alpha}^{\pm}$. We have to split up into two cases according to the parity of $d_{\alpha}^{+}\,d_{\alpha}^{-}$. If it is even, then exactly one of $d_{\alpha}^{+}$ and $d_{\alpha}^{-}$ is even and we find that $c_{\psi} = 2 n_{\alpha}^{+}n_{\alpha}^{-}= n_{\psi}^{+} n_{\psi}^{-}$. Otherwise, if it is odd, $c_{\psi} = 2 n_{\alpha}^{+}n_{\alpha}^{-}= \tfrac{1}{2} n_{\psi}^{+} n_{\psi}^{-}$. It remains to note that $d_{\alpha}^{+}\,d_{\alpha}^{-}$ is odd if and only if $n_{\psi}^{+}$ and $n_{\psi}^{-}$ are both even.

 Finally, we can treat the case when $E_0(\RR)$ has one and $E(\RR)$ has two connected components by a similar case-by-case treatment. Alternatively one can just apply the above to the dual of $\beta$.
\end{proof}

There are other ways to find $c_{\alpha}$ for an isogeny $\alpha\colon E\to E'$. For instance, the expansion of $\alpha$ using the formal groups for $E$ and $E'$ will have $c_{\alpha}$ as the leading coefficient. Also there is the useful formula $c_{\alpha}^2 = \deg(\alpha)\, c_{\infty}(E)\, \Omega^{+}_E\,\Omega_E^{-} / (c_\infty(E') \,\Omega^+_{E'}\, \Omega^-_{E} )$. The advantage of the formula in Proposition~\ref{manin_prop} is that all terms can be read off $E$ and $E'$ without reference to $\psi$ any more.

  As an example we add here the case of the isogeny class 27a. There are four curves in this class and they are linked by the following $3$-isogenies
  \begin{equation*}
    \xymatrix@R-7ex{&& \text{27a3}\ar[rd]\ar[ld] & \\
    & \text{27a1}\ar[dl] &  & \text{27a4}\\
    \text{27a2} }
  \end{equation*}
  where the direction of the arrow indicates the isogeny $\alpha$ for which $c_{\alpha}=1$. In other words, the inclusion of the N\'eron lattices is in the opposite direction. The curve 27a1 is the strong Weil curve, while 27a3 is the minimal curve in the sense of~\cite{stevens}. The three curves on the right have each exactly $3$ points in $E(\QQ)$ and they lie in the kernel of the isogeny to the curve on their left. The Manin constants are equal to $1$ for 27a2 and 27a1 and they are equal to $3$ for the two curves 27a3 and 27a4.

\subsection{Outstanding issues}\label{outstanding_subsec}

  There are two outstanding issues. First, what happens if $c_0\neq 1$ and secondly how do we find the strong Weil curve in the isogeny class.

  Suppose that the Manin constant $c_0$ were larger than $1$. If we knew the value of $c_0$ we could simply multiply the bounds $t_0$ and $t_0^{\pm}$ by $c_0$, too. However, it is then likely that we would at some point find a modular symbol where $c_0$ appears as a factor of the denominator. When rounding our numerical approximation, we would find a large error. If this happens, we could verify that $c_0\neq 1$ and announce the exceptional news to the world. Therefore we do not really have to worry about this assumption in practice.

  For all isogeny classes in Cremona's tables~\cite{cremona} it has been verified that $c_0=1$ when the table was created. For a few curves this is slightly more complicated and the issue is well explained in the appendix of~\cite{manin_constant}.

  The second issue is related to the first. Even for the curves in the tables, it is not always possible to say with certainty which curve in the isogeny class is the strong Weil curve. This arises because the computation in creating the table is done mostly with $+$-modular symbols only. At worst, we are off by a lattice of index $2$.

  Finally, suppose the curve lies outside the range of the table. We can still determine the isogeny class of the curve fairly quickly. However we have no means of knowing which curve is the strong Weil curve. To be on the safe side, we have to assume that it is one of the curves with maximal lattice. In practice it is very often on the contrary the minimal curve that is the strong Weil curve, but we have no way of showing this for our curve. If we are really unlucky, we even picked the wrong curve among the maximal curves; hence we should really work with the lattice generated by all N\'eron lattices in the isogeny class.


\section{Numerical integration}\label{num_sec}

  Let $f$ be the newform associated to the isogeny class of the elliptic curve $E$. Let $\ve>0$.  In this section, we consider the finite sum that approximates the integral of $2\pi i f(z)dz$ from $i\infty$ to a point $\tau$ in the upper half plane. We prove bounds on the number of terms and the bit precision to work with in order to determine the integral with an error of at most $\ve$.

  Generalising the definition of $\lambda(r)$, we will consider
  \begin{equation*}
    \lambda(\tau) = 2\pi i \int_{i\infty}^{\tau} f(z) dz
  \end{equation*}
  for any point $\tau = x + yi$ in the upper half plane. As $y>0$, we can express it as the evaluation of a power series in $q = e^{2\pi i \tau}$, namely
  \begin{equation}\label{sum_lambda_eq}
    \lambda(\tau) = \sum_{n=1}^{\infty} \frac{a_n}{n} q^n = \sum_{n=1}^\infty \frac{a_n}{n} \exp(-2\pi ny + 2\pi n x i).
  \end{equation}
  We will approximate this sum by its finite partial sum for $n\leq T$ for a bound $T$. It is the value of $y$ that determines how quickly the sum will converge and so how large $T$ should be. In Section~\ref{partial_sec}, we will be interested  in the following partial sums: for any $m>1$ and $0\leq j <m$ and $y>0$, we define
  \begin{equation}\label{sum_kappa_eq}
    \kappa_{j,m} (y) = \sum_{\substack{n\geq 1\\ n\equiv j \bmod{m}}} \frac{a_n}{n} \exp(-2\pi n y) \in\RR.
  \end{equation}

\subsection{Truncation}\label{truncation_subsec}

  We now proceed to determine how many terms in the sums in~\eqref{sum_lambda_eq} and~\eqref{sum_kappa_eq} we have to add to be guaranteed a value that differs from the infinite sum by less than a given error $\ve$. Afterwards we will decide with what level of precision we have to do the numerical computations so that the error due to precision loss will be smaller than a given bound $\ve'$. Recall that we have determined the value of $\ve+\ve'$ in Section~\ref{denominator_sec}.

  Define the following function for $y>0$ and $\ve>0$.
  \begin{equation}\label{number_of_terms_eq}
    T(y,\ve) = \frac{-\log(2 \pi y \ve)}{2\pi y}
  \end{equation}
  which is, for a fixed $\ve$, a function that grows like a constant multiple of
  $\tfrac{1}{y}\log(\tfrac{1}{y})$ as $y\to 0$.

  \begin{lem}\label{truncate_lem}
    Let $\tau$ be an element of the upper half plane with $y =\im(\tau)$ and let $\ve>0$. If
    $T> T(y,\ve)$ then we have
    \begin{equation*}
      \Biggl\vert \lambda(\tau) - \sum_{n=1}^{T} \frac{a_n}{n} \exp(2\pi i n \tau)\Biggr \vert < \ve\,.
    \end{equation*}
  \end{lem}
  \begin{proof}
    Write $\tau = x+yi $. Now we use that $\vert a_n\vert \leq n$ as proven in Lemma~2.9 in~\cite{stein_bsd}. The difference to bound is
    \begin{equation*}
      \Biggl\vert \sum_{n > T} \frac{a_n}{n} \exp\Bigl(2\pi i n (x+i\, y)\Bigr) \Biggr\vert
      \leq \sum_{n>T} \frac{\vert a_n\vert}{n} \exp(-2\pi  n y) \leq \sum_{n>T}\exp(-2\pi  n y) = \frac{e^{-2\pi (T+1)y}}{1-e^{-2\pi y}} =\frac{e^{-2\pi T y}}{e^{2\pi y}-1}.
    \end{equation*}
    Now the condition on $T$ implies that
    \begin{equation*}
      \frac{e^{-2\pi Ty}}{e^{2\pi y}-1} < \frac{2\pi y \ve}{e^{2\pi y} - 1} < \ve.\qedhere
    \end{equation*}
  \end{proof}
  In this proof, we have used the inequality $\vert a_n \vert \leq n$.  In fact, we even know that $\vert a_n\vert \leq \sigma_0(n)\,\sqrt{n}$ where $\sigma_0(n)$ is the number of positive divisors of $n$. However even this asymptotically sharper inequality will not lead to a substantially better theoretical bound on the number of terms.

  Nonetheless, in practice we use the following estimates. First we have the trivial bound $\sigma_0(n)\leq 2\sqrt{n}$. Moreover for every $2>\varsigma > 0$ the equality $\sigma_0(n) < \varsigma \cdot\sqrt{n}$ holds for all $n>B(\varsigma)$ for some $B(\varsigma)$. Here are a few values of this bound used in the implementation:
  \begin{center}
      \begin{tabular}{r|*{7}{c}}
          $\varsigma$    & $1$ & $2/3$ & $1/2$ & $1/3$ & $1/4$ & $1/5$ & $1/6$ \\
          $B(\varsigma)$ & $1260$ & $10080$ & $55440$ & $277200$ & $831600$ & $2162160$ & $4324320$
      \end{tabular}\qquad\ \refstepcounter{equation}\label{Bvarsigma_eq}$(\theequation)$
  \end{center}

  With the same method as in Lemma~\ref{truncate_lem} one proves the bound on the approximation for the partial sum $\kappa_{j,m}(y)$.
  When we will compare the methods it will be clear that the corresponding error bound that we ask for is $\ve/m$.
  \begin{lem}\label{Tprime_lem}
    Let $y>0$, $m>1$, $0\leq j <m$ and $\varepsilon > 0$. If $T > T(y,\ve)+m$, then
    \begin{equation*}
      \Biggl\vert \kappa_{j,m}( y) - \sum_{\substack{n\equiv j \bmod{m}\\ 1\leq n \leq T}} \frac{a_n}{n} \exp(-2\pi n y) \Biggr\vert < \frac{ \varepsilon}{m}.
    \end{equation*}
  \end{lem}
  We have seen that the value of $\tfrac{1}{y}$ is an important measure of how difficult it will be to approximate the integral. This motivates the following definition.
  \begin{defn}
   We call the value of $y$ the \textit{speed} of the evaluation of $\lambda(x+yi)$.
  \end{defn}
  The larger the speed the faster we can compute $\lambda(\tau)$.

  Of course, since the sums are alternating in average (because the $a_p$ for primes $p$ follow the Sato-Tate distribution), they actually converge much faster. In~\cite{goldfeld}, Goldfeld suggests that it is probable that the computation complexity is polynomial in $N$; in other words that the bound for $T$ could behave like a power of $\log(\tfrac{1}{y})$. However this is still far beyond current knowledge. Even an unproven effective version of the Sato-Tate distribution does not seem to help here.

\subsection{Implementation}
  For implementing these finite sums we use Horner's rule. Here is the algorithm to evaluate an approximation to $\lambda(\tau)$. We are given $\tau$ in the upper half plane and an bound $\ve$ on the allowed error.

  \begin{algo}{Algorithm: Numerical approximation to $\lambda(\tau)$.}
    \item[Initialisation] Set $s\leftarrow 0$ and $n\leftarrow \lceil T(y,\ve)\rceil $ and compute $q \leftarrow \exp( 2 \pi i \tau)$.
    \item[Loop] While $n$ is positive, replace $s\leftarrow s\cdot q + \tfrac{a_n}{n}$ and decrease $n$ by one.
    \item[End] Return $s\cdot q$ as a good approximation to $\lambda(\tau)$.
  \end{algo}
  The same idea can be used to compute an approximation to the partial sum $\kappa_{j,m}(y)$ for all $j$ simultaneous. We are given $m$ and $y$ and the allowed error $\ve/m$.
  \begin{algo}{Algorithm: Simultaneous numerical approximation to $\kappa_{j,m}(y)$.}
    \item[Initialise] Set $v_j\leftarrow 0$ for all $0\leq j< m$. Compute $q \leftarrow \exp(-2\pi y)$ and $q' \leftarrow \exp(-2\pi m y)$. Set to start $n\leftarrow \lceil T(y,\ve)\rceil $.
    \item[Loop] As long as $n$ is positive, replace $v_j\leftarrow v_j\cdot q' +\tfrac{a_n}{n}$, where $j\equiv n \pmod{m}$ and then decrease $n$ by $1$.
    \item[End]  At the end the value $v_j\cdot q^j$ for $1\leq j <m$ and $v_0 \cdot q'$ are good approximations to $\kappa_{j,m}(y)$.
  \end{algo}

\subsection{Precision}\label{prec_subsec}
  We wish to determine with how many bits $b$ of precision we have to work with to make sure that the error in the above algorithm is smaller than a given error $\ve'$. In practice this error will be chosen to be a tiny fraction of the error $\ve$ that we allowed for finding the above bound $T$.

  \begin{lem}\label{prec_lem}
    Let $1>\varepsilon'>0$ and $\tau$ in the upper half plane. Let $T$ be the number of terms evaluated in the sum to approximate $\lambda(\tau)$. If
    \begin{equation*}
      2^{-b} < \frac{\varepsilon'}{2 T(T+\varepsilon')}
    \end{equation*}
    then the numerical value computed differs from the actual sum $\sum_{n=1}^T \tfrac{a_n}{n} q^n$ by less than $\ve'$ in absolute value.
  \end{lem}
  \begin{proof}
    We may suppose that the value of $q$ can be pre-computed to $b$ bits of precision. We use the absolute error estimate on the Horner's rule given on page~105 of~\cite{higham}. If we write $\delta = 2^{-b}$, then the absolute error is smaller than
    \begin{equation*}
      \frac{2T\delta}{1-2T\delta} \cdot \sum_{n=1}^T \Bigl\vert \frac{a_n}{n}\Bigr \vert \cdot \vert e^{2\pi i \tau}\vert ^n \leq \frac{2T\delta}{1-2T\delta} \cdot \sum_{n=1}^T e^{-2\pi ny} \leq \frac{2T\delta}{1-2T\delta} \cdot T,
    \end{equation*}
    where we used again that $\vert a_n\vert \leq n$. It is now easy to see that the given inequality on $\delta$ in the lemma implies that the above right hand side is smaller than $\varepsilon'$
  \end{proof}
  For the approximation of $\kappa_{j,m}(y)$ to have an error smaller than $\ve'/m$, we have to impose the bound
  \begin{equation*}
    2^{-b} < \frac{\varepsilon'\, m}{2T'(T'+ \varepsilon')}
  \end{equation*}
  where $T' = T(y,\varepsilon)+m$ is the upper limit of the finite sums in Lemma~\ref{Tprime_lem}.

  Later, it will be clear later that, in view of Lemma~\ref{prec_lem}, we may neglect the issue of memory usage because the floating point numbers will take up approximatively as many bits as the conductor or the coefficients of $E$ take up.

  Within the range of interesting examples, the standard double precision of 53 bits is often sufficient. For example, the period $\Omega^{+}$ of the curve 100002a1 is approximatively equal to $1.125$. If we set $\varepsilon=0.278427$ and $\varepsilon'=0.002812$, then we are allowed to sum up $T=3558923$ terms using $53$ bits precision, which would allow for $\tfrac{1}{y}$ to be as large as $1627105$. From the results in the following sections one can deduce that this allows to evaluate all Manin symbols using standard double precision.

  Instead, for a curve like $E\colon y^2 = x^3+101x+103$ of conductor $35261176$, the evaluation of $\lambda(\tfrac{1}{107})$ will require precision above $53$ bits to obtain provable results.


\section{Computation of unitary symbols}\label{unitary_sec}

  In this section we assume that $r$ is a unitary cusp. It is equivalent to the definition given at the start of Section~\ref{denominator_sec} to ask that the cusp $r$ on $X_0(N)$ is in the orbit of $i\infty$ under the group of Atkin-Lehner involutions. This section explains how to compute $\lambda(r)$ under this assumption. In Section~\ref{num_sec}, we explained how to compute integrals from a point $\tau$ within the upper half plane to the cusp $i \infty$. Now, we wish to explain how one integrates paths from $\tau$ to another cusp $r\in\QQ$. The idea to use the Atkin-Lehner involution to bring $r$ to $\infty$ is already presented in~\cite{goldfeld}.

\subsection{Moving unitary cusps with Atkin-Lehner involutions}\label{moving_subsec}

  By assumption, $r=\tfrac{a}{m}$ is unitary. Recall that we denote by $M$ be the greatest common divisor of $m$ and the conductor $N$. Further we write $N = Q\cdot M$ with $Q$ and $M$ coprime. Then the greatest common divisor of $Qa$ and $m$ is $1$ and hence we find integers $u$ and $v$ such that $Qau + mv = 1$. We define
  \begin{equation*}
    W_r = \begin{pmatrix} Qu & v \\ -Qm & Qa \end{pmatrix}
  \end{equation*}
  which is of determinant $Q$ and sends $r$ to $i\infty$ under the action of $\GL_2(\QQ)$ on the completed upper half plane. Since $Qm$ is divisible by $N$, the matrix $W_r$ induces an Atkin-Lehner involution on $X_0(N)$. Since $f$ is a newform it is also an eigenfunction for $W_r$. We have $f\vert_{W_r} = \epsilon_Q \cdot f$ for $\epsilon_Q\in\{\pm 1\}$. In fact, $\epsilon_Q$ is easy to compute as it is just the product of the local root numbers for $\ell\mid Q$; and for a product of semistable primes, we have $\epsilon_Q = -a_Q$. We get
  \begin{equation}\label{int_r_tau_eq}
    2\pi i \int_{\tau}^r f(z)dz = \epsilon_Q \cdot 2\pi i \int_{\tau}^r f\vert_{W_r} (z) dz = \epsilon_Q  \cdot 2\pi i \int_{W_r(\tau)}^{i \infty} f(z) dz = -\epsilon_Q \cdot \lambda\bigl(W_r(\tau)\bigr)
  \end{equation}
  which can be evaluated with the previously described numerical method. Note that the speed of this evaluation is equal to
  \begin{equation}\label{wr_speed_eq}
    \im\bigl(W_r(\tau)\bigr) = \frac{Q\cdot \im(\tau)}{\vert -Qm\tau +Qa\vert^2} = \frac{\im(\tau)}{Q \cdot m^2 \cdot \vert r-\tau\vert ^2}.
  \end{equation}

\subsection{Splitting up the path from \texorpdfstring{$i\infty$}{ioo} to \texorpdfstring{$r$}{r} }\label{ioo_to_cusp_subsec}

  We wish to compute $\lambda(r)$ by splitting up the path of integration from $r$ to $i\infty$ at a certain $\tau$ in the upper half plane.
  Using~\eqref{int_r_tau_eq}, we find, for any such $\tau$ and any unitary cusp $r$,
  \begin{equation}\label{split_r_ioo_eq}
    \lambda(r) = 2\pi i \Bigl(\int_{i\infty}^{\tau} + \int_{\tau}^r \Bigr) f(z)dz
               =  \lambda(\tau)- \epsilon_Q\cdot \lambda\bigl(W_r(\tau)\bigr).
  \end{equation}
  These two values of $\lambda$ can be evaluated using the numerical method. We are now looking for the choice of $\tau$ such that the computation is fastest. The following lemma will show that this is achieved when the speed of computing $\lambda(\tau)$ is equal to the speed of computing $\lambda\bigl (W_r(\tau)\bigr)$ and they are maximal. From~\eqref{wr_speed_eq}, the first condition is equivalent to the equation
  \begin{equation*}
    \im(\tau) = \frac{\im(\tau)}{Q \cdot m^2 \cdot \vert r-\tau\vert ^2}.
  \end{equation*}
  So we are looking for the $\tau= x +yi$ with maximal $y$ such that $\vert r - \tau\vert =1/\bigl(\sqrt{Q}\,m\bigr)$. This is obtained for
  \begin{equation*}
    \tau = r + \frac{1}{m\sqrt{Q}} \, i.
  \end{equation*}
  It is not difficult to see that
  \begin{equation}\label{Wrtau_eq}
    W_r \Bigl(\frac{a}{m} + \frac{1}{m\sqrt{Q}}\,i\Bigr) = -\frac{u}{m} + \frac{1}{m\sqrt{Q}}\,i
  \end{equation}
  where $u$ is an inverse of $Q\cdot a$ modulo $m$.
  We still have to justify the claim that our choice of $\tau$ is optimal.
  \begin{lem}
    For a fixed curve $E$, a fixed unitary cusp $r$, and a fixed error $\ve >0$, the minimum of
    \begin{equation*}
      T\Bigl(y, \frac{\ve}{2}\Bigr) + T\Bigl( \frac{y}{Q m^2 \vert r - \tau\vert^2}, \frac{\ve}{2}\Bigr)
    \end{equation*}
    is attained when $\tau = r + \frac{1}{m\sqrt{Q}} \,i$.
  \end{lem}
  \begin{proof}
    Write as before $\tau = x + y\,i$. Since $t(y) := T(y,\ve/2)$ is decreasing in $y>0$, the best choice for $x$ must occur when $x = r$. The function to minimise simplifies then to $t(y) + t(1/(Cy))$ with $C= Q\,m^2$. Taking the derivatives with respect to $y$, we see that at the minimum, we must have
    \begin{equation*}
      y \cdot t'(y) = \frac{1}{Cy} \cdot t'\Bigl(\frac{1}{Cy}\Bigr).
    \end{equation*}
    Now from the definition we see that $t$ satisfies the differential equation
    \begin{equation*}
      y\cdot  t'(y) = - t(y) - \frac{1}{2\pi y}
    \end{equation*}
    and hence, since $t$ is decreasing, $y\mapsto y\cdot t'(y)$ is increasing. Hence there is only one minimum, namely when $y = \frac{1}{Cy}$.
  \end{proof}


\subsection{Integrals from cusp to cusp}\label{cusp_to_cusp_subsec}

  Let $r=\tfrac{a}{m}$ and $r'=\tfrac{a'}{m'}$ be two unitary cusps of widths $Q$ and $Q'$ respectively. Our aim is to compute
  \begin{equation*}
    \lambda\bigl( \{r'\to r\} \bigr) = 2\pi i \int_{r'}^r f(z) dz = \lambda ( r ) - \lambda(r'),
  \end{equation*}
  where the integration follows any path from $r'$ to $r$ in the upper half plane.
  One way to do so, indicated by the last expression above, is to integrate from $r'$ to $i \infty$ and then subtract the integration from $r$ to $i\infty$ using the method explained above. We call this the \textit{indirect} way.

  Instead, the \textit{direct} way splits up the integration path from $r'$ to $r$ into two pieces: First find a good $\tau$ in the upper half plane. Then use the Atkin-Lehner involution $W_{r'}$ to move the path from $r'$ to $\tau$ to a path from $i\infty$ to $W_{r'}(\tau)$. Similarly use $W_r$ to move the second piece to a path from $W_r(\tau)$ to $i\infty$. As before, on these two paths we can use the methods from the previous section. We get
  \begin{equation*}
    \lambda\bigl( \{ r'\to  r\} \bigr) = 2\pi i \biggl(\int_{r'}^{\tau}+\int_{\tau}^r \biggr) f(z) dz =
    \epsilon_{Q'} \cdot \lambda\bigl( W_{r'}(\tau)\bigr ) - \epsilon_{Q}\cdot \lambda\bigl(W_{r}(\tau) \bigr)\,.
  \end{equation*}
  We expect again the best choice for $\tau=x+i\,y$ to be such that the speeds of both integrals are equal and they are maximal. We get the equation
  \begin{equation}\label{compare_speeds_eq}
    \frac{y}{m^2 Q \cdot \vert \tau -r \vert^2} = \frac{y}{m'^2 Q' \cdot \vert \tau -r' \vert^2}.
  \end{equation}
  If we denote
  \begin{equation*}
    c = \frac{m'}{m}\cdot \sqrt{\frac{Q'}{Q}} >0
  \end{equation*}
  then the above equation~\eqref{compare_speeds_eq} becomes
  \begin{equation*}
    \vert \tau -r\vert = c \cdot \vert \tau - r'\vert.
  \end{equation*}
  The set of all complex numbers satisfying this equation forms a circle around either $r$ or $r'$. More precisely, if $c>1$ then it is a circle with $r$ in its interior and $r'$ in the exterior. Conversely if $c<1$ then it is $r'$ that lies in the interior and $r$ in the exterior. Finally if $c=1$, we deal with a vertical line bisecting the segment from $r$ to $r'$. Write $A=\vert \tau - r\vert$ and $A'=\vert \tau -r'\vert$.
  \begin{center}
    \includegraphics[width=14cm]{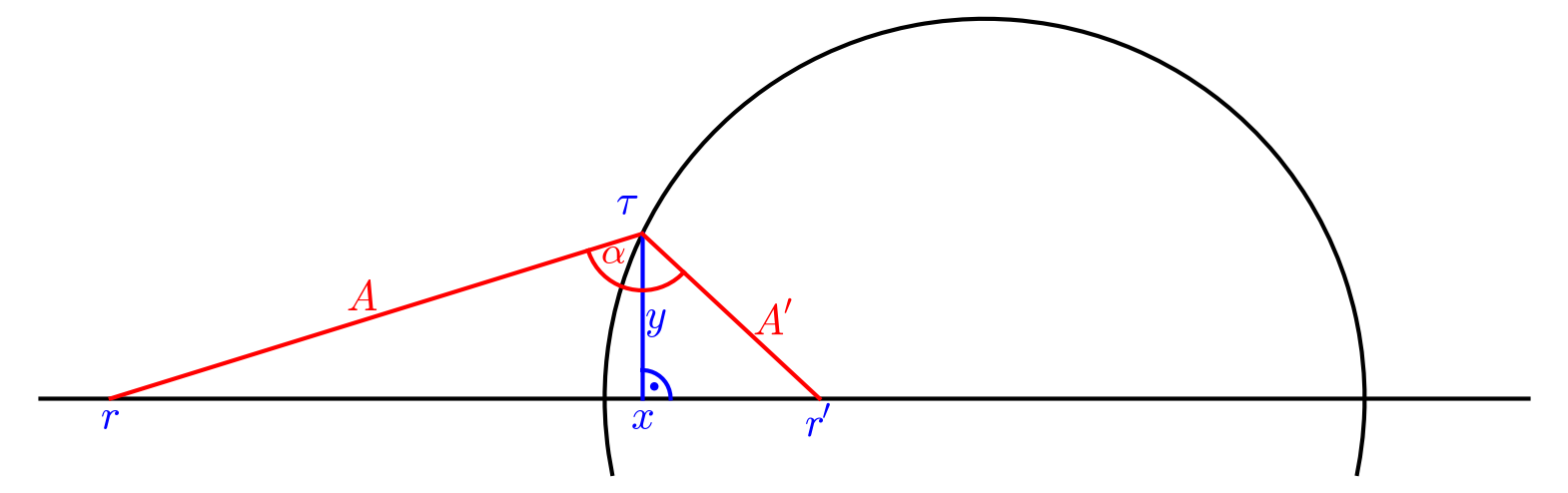}
  \end{center}
  Our aim now is to maximise the function in~\eqref{compare_speeds_eq}, which is the same as to maximise $\frac{y}{A^2}$, on this circle $A=c\cdot A'$. We have
  \begin{equation*}
    \frac{y}{A^2} = \frac{y}{A\cdot A'\cdot c} = \frac{\sin(\alpha)}{\vert r' - r \vert\cdot c}
  \end{equation*}
  where $\alpha$ is the acute angle between the segments from $\tau$ to $r$ and $r'$ respectively. This is maximal when $\alpha=\pi/2$. So $\tau$ is the intersection of the circle $A=c\cdot A'$ with the circle centred on the real axis and passing through $r'$ and $r$. It is now easy to compute that
  \begin{align*}
    y &= \frac{c}{c^2+1} \cdot \vert r'-r\vert = \frac{\sqrt{Q\, Q'}}{m^2\,Q+m'^2\,Q'}\cdot\vert a\,m'-a'\,m\vert \\
    x &= \frac{c^2 r'+r}{c^2+1} =\frac{a\,m\,Q+a'\,m'\,Q'}{m^2\,Q+m'^2\,Q'}
  \end{align*}
  The maximum value for the speed in~\eqref{compare_speeds_eq} is
  \begin{equation*}
    \frac{1}{m\,m'\sqrt{Q\,Q'}\cdot \vert r-r'\vert} = \frac{1}{\sqrt{QQ'} \cdot\vert am'-a'm\vert}.
  \end{equation*}
  Furthermore, we find
  \begin{align*}
    W_r(\tau) &= \frac{1}{Q}\frac{Qa'u+m'v}{am'-a'm} + \frac{i}{\sqrt{QQ'} \cdot\vert am'-a'm\vert}\qquad\text{and} \\
    W_{r'}(\tau) &= \frac{1}{Q'}\frac{Q'au'+mv'}{a'm-am'} + \frac{i}{\sqrt{QQ'} \cdot\vert am'-a'm\vert}.
  \end{align*}
  We could not spot any general rule to distinguish the cases when the direct or the indirect method is faster. In practice it is easy to test before starting to sum. For the curve $E$=5077a1 and $\ve =0.001$, the direct method is faster for $(r,r')= (0,\tfrac{70}{5077})$, but slower for $(r,r') = (\tfrac{123}{456},\tfrac{789}{5077})$.


\section{Computation of non-unitary symbols}\label{non_unitary_sec}

   If $N$ is not square-free then there are modular symbols that we do not know how to compute with the above methods. In the section, we analyse how to compute $\lambda(r)$ when $r$ is non-unitary. We cannot move the cusp to $i\infty$ using an Atkin-Lehner involution. If the elliptic curve admits a quadratic twist $E^{\dagger}$ whose conductor is square-free, then it is best to use the formula for twisting modular symbols, see Section~\ref{twists_subsec}. But this is not always possible.

  There is one special case when we can transform a non-unitary symbol to a unitary one: Suppose $4\mid N$ and $r=\tfrac{a}{2m}$ with odd $a$ and $m$. Then the action of the Hecke operator $T_2$ yields the equality
  \begin{equation*}
    \lambda\Bigl(\frac{r'}{2}\Bigr) + \lambda\Bigl(\frac{r'-1}{2}\Bigr)=0
  \end{equation*}
  because $a_2=0$. For $r'/2=r$, we get $\lambda(r) = -\lambda\bigl(\tfrac{a-m}{2m}\bigr)$. The latter is now at a cusp with an odd denominator and has a chance of being a unitary cusp. This little trick only works for $4\mid N$ not any other square dividing $N$.

  In general, however, we know no better method than to rewrite $\lambda(r)$ as the sum of so-called transportable symbols via the use of a Hecke operator. We start by explaining what transportable symbols are and how they can be computed.

\subsection{Transportable modular symbols}\label{transportable_subsec}

 \begin{defn}
   We will call $\lambda(\{r'\to r\})$ a \textit{transportable} modular symbol if the
   two rational numbers $r$ and $r'$ are $\Gamma_0(N)$-equivalent.
 \end{defn}
  This is a more restrictive definition of this term than in~\cite{stein_verrill} where they allow also sums of transportable symbols in the more general setting of higher weight modular forms.

  Let $\lambda(\{r'\to r\})$ be a transportable modular symbol. We may compute it by transporting the path: if $\gamma\in\Gamma_0(N)$ is such that $r'=\gamma(r)$, then
  \begin{equation}\label{transportable_eq}
    \lambda\bigl(\{r'\to r\}\bigr) = 2\pi i\int_{\gamma(r)}^{r} f(z) dz = 2 \pi i \int_{\gamma(\tau)}^{\tau} f(z) dz = \lambda(\tau) - \lambda(\gamma(\tau))
  \end{equation}
  for any $\tau$ in the upper half plane.

   Write $\gamma = \sm{a}{b}{c}{d}$. Note first that if $\gamma$ is not hyperbolic, i.e., if $\vert a+d\vert \leq 2$, then there is a point $\tau$ in the upper half plane or among the cusps with $\gamma(\tau)= \tau$ and thus $\lambda\bigl(\{r'\to r\}\bigr)=0$. Hence we may assume that $\gamma$ is hyperbolic.

  Let us now find the best choice of $\tau = x + yi$ in the upper half plane. It will be such that the speeds of summing up $\lambda(\tau)$ and $\lambda(\gamma(\tau))$ are equal and as large as possible. This implies that
  \begin{equation*}
    (cx +d)^2 + c^2 y^2 = 1.
  \end{equation*}
  We want to maximise $y$ under this restriction, so obviously the best choice is $y=\tfrac{1}{\vert c\vert}$ and $x=-\tfrac{d}{c}$ and the speed will be $1/\vert c\vert$. See Algorithm 10.6 in~\cite{stein}. Since $c\in N\ZZ$, the speed will be smaller than $\frac{1}{N}$, which is often quite worse than the previous methods.

  Given two $\Gamma_0(N)$-equivalent cusps $r$ and $r'$, we should try to find the matrix $\gamma\in\Gamma_0(N)$ with $\gamma(r)=r'$ in such a way as to make its lower left entry $c$ as small as possible in absolute value. Proposition 2.2.3 in~\cite{cremona} gives an algorithm to construct such a matrix, but reading carefully the proof one sees that it actually gives the construction of all possible $\gamma$. We repeat it here in our notations for the convenience of the reader.

  Write $r=\tfrac{e}{m}$ and $r'=\tfrac{e'}{m'}$ in reduced fractions. Using the euclidean algorithm, we can find matrices $\delta =\sm{e}{u}{m}{v}$ and $\delta' = \sm{e'}{u'}{m'}{v'}$ in $\SL_2(\ZZ)$ such that $\delta( i \infty) = r$ and $\delta'(i \infty) = r'$. We have that $\gamma_0=\delta'\cdot\delta^{-1} = \sm{a_0}{b_0}{c_0}{d_0}$ sends $r$ to $r'$. We can obtain all such matrices as
  \begin{equation*}
    \gamma = \delta'\cdot \begin{pmatrix} 1 & t \\ 0 & 1 \end{pmatrix} \cdot \delta^{-1}= \begin{pmatrix} a &  b \\  c &  d \end{pmatrix}
  \end{equation*}
  for some $t\in\ZZ$, because $\delta'^{-1}\gamma\delta$ stabilises the cusp $i \infty$. (Alternatively one can view this as the indeterminacy of $u$ and $v$ in the B\'ezout equation $e v - m u = 1$ modulo $e$ and $m$ respectively.)  The equation
  \begin{equation*}
    0\equiv c = c_0 + t\, m\, m' \pmod{N}
  \end{equation*}
  is solvable in $t$ if and only if $r$ and $r'$ are $\Gamma_0(N)$-equivalent. The choice of $c$ is unique up to multiples of $\lcm(mm',N)$. So we can take $t$ such that $c$ is the least residue modulo $\lcm(mm',N)$; hence we have that the speed will be at least $\frac{2}{\lcm(mm',N)}$.

  Now so far, we have considered to transport the path close to the cusp $i\infty$. However, we could also choose another unitary cusp $r_0$ of width $Q$ instead. We compute
  \begin{equation*}
    \lambda\bigl(\{r'\to r\}\bigr) = 2 \pi i \Biggl(\int_{r_0}^{\tau} -\int_{r_0}^{\gamma(\tau)}\Biggr) f(z) dz = \epsilon_Q \Bigl( \lambda\bigl(W_{r_0}(\tau)\bigr) - \lambda\bigl(W_{r_0}(\gamma(\tau))\bigr)\Bigr)
  \end{equation*}
  by~\eqref{int_r_tau_eq}. Renaming $W_{r_0}(\tau)$ as $\tau$ and writing $\gamma_{r_0} = W_{r_0}\cdot \gamma \cdot W_{r_0}^{-1}\in \Gamma_0(N)$, this is equal to
  \begin{equation*}
    \lambda\bigl(\{r'\to r\}\bigr) =  \epsilon_Q \Bigl( \lambda(\tau) - \lambda\bigl(\gamma_{r_0}(\tau)\bigr)\Bigr)
  \end{equation*}
  for all $\tau$ in the upper half plane. The best $\gamma$ is obtained when the lower left entry of $\gamma_{r_0}$ is minimal. Again, this entry is divisible by $N$ and we expect a rather low speed.

  For example, we can take $r_0 = 0$ of width $N$. Then $\gamma_{r_0} =\sm{d}{-c/N}{-bN}{a}$ and so we are now looking for $\gamma$ such that $\vert b\vert$ is minimal. As before $b = b_0 + t\, e\, e'$ and we are looking for the least residue of $b_0$ modulo $ee'$. It may be that the resulting computation is faster with $r_0$ than with $i\infty$. It seems difficult to find the best choice of the unitary cusp $r_0$ in general.

  Finally, we could also transport the path in such a way as to have $\gamma(\tau)$ close to $i\infty$ and $\tau$ close to another unitary cusp $r_0$. For instance if $r_0=0$, this would give a speed of $1/\bigl(\sqrt{N}\vert d\vert\bigr) $ and we would have to minimise $\vert d\vert$. However, this time it will also involve the computation of the integral from $i\infty$ to $r_0$.

\subsection{Hecke operators to get transportable paths}

  Let $r=\tfrac{a}{m}$ be a non-unitary cusp. Set $M$ to be the greatest common divisor of $m$ and $N$. Further put $d$ equal to the greatest common divisor of $M$ and $Q= \frac{N}{M}$. The previous methods explain how to compute $\lambda(r)$ only in the case that $d=1$. In this section, we will suppose $d>1$.

  First, for any integer $n$ coprime to $N$, we have the action of the Hecke operator, which gives us
  \begin{equation*}
    a_n \cdot \lambda(r) = \sum_{k\mid n} \sum_{u=0}^{k-1} \lambda \Bigl(\frac{nr+uk}{k^2} \Bigr).
  \end{equation*}
  The cusp $\frac{nr+uk}{k^2}$ is $\Gamma_0(N)$-equivalent to $r$ if and only if $n\cdot k^{-2}$ is congruent to $1$ modulo $d$. This implies that $n\equiv 1 \pmod{d}$ and that $k^2\equiv 1 \pmod{d}$ for all divisors $k\mid n$. If $n$ is not a prime or a square of a prime, then the smallest prime divisor of $n$ will provide a smaller choice for $n$.

  Let $\ell$ be a prime congruent to $1$ modulo $d$. If $\ell$ does not divide $N$  then
  \begin{equation*}
    \bigl(a_{\ell}-\ell - 1 \bigr) \lambda(r) =
    \lambda\bigl(\{\ell r\to r\}\bigr) +\sum_{u=0}^{\ell-1} \lambda\biggl(\Bigl\{\frac{r+u}{\ell}\to r \Bigr\}\biggr).
  \end{equation*}
  The right hand side is now a sum of $\ell+1$ transportable symbols. The integer $a_{\ell}-\ell-1=-N_{\ell}$ is non-zero since $N_{\ell}$ is the number of points on the reduction of $E$ to $\FF_{\ell}$.
  If $\ell$ divides $N$, then we get
   \begin{equation*}
    \bigl(a_{\ell}-\ell  \bigr) \lambda(r) =
    \sum_{u=0}^{\ell-1} \lambda\biggl(\Bigl\{\frac{r+u}{\ell}\to r \Bigr\}\biggr)
  \end{equation*}
  instead. This time $\vert a_{\ell}\vert \leq 1$.

  The other option is to take a prime $\ell$ such that $\ell^2\equiv 1\pmod{d}$. For instance, let $\ell\equiv -1\pmod{d}$. Then we have the following formula
  \begin{equation*}
    \bigl(a_{\ell^2} -\ell^2 -\ell- 1\bigr) \lambda(r) = \lambda\Bigl(\{\ell^2 r\to r\}\Bigr)
    + \sum_{u=0}^{\ell-1} \lambda\biggl(\Bigl\{ r+\frac{u}{\ell}\to r\Bigr\}\biggr)
    + \sum_{v=0}^{\ell^2-1} \lambda\biggl(\Bigl\{\frac{r+v}{\ell^2} \to r\Bigr\}\biggr)
  \end{equation*}
  which expresses a non-zero multiple of $\lambda(r)$ as a sum of transportable symbols. If $\ell$ is the smallest prime congruent to $1$ modulo $d$ and $\ell'\not\equiv 1\pmod{d}$ is the smallest prime such that $\ell'^2\equiv 1\pmod{d}$, then the above formula for $\ell'$ will have $\ell'^2+\ell'+1$ terms, which may be smaller than the $\ell+1$ terms in the corresponding sum for $\ell$. Although not frequent, there are cases when this is useful. For instance if $d=6441$, we have $\ell = 231877$ and $\ell' = 227$.

  It is hard to estimate what the complexity of this method is. It is certainly significantly slower than the computation of unitary cusps, but it is still useful when the conductor is not too large. In the most frequent applications, like for the computation of $p$-adic $L$-series, this is not important, as we will be mainly interested in unitary symbols. Note however that the following section shows that even the computation of unitary symbols for large denominators may encounter the computations explained here.


\section{Manin's trick using continued fractions}\label{manins_trick_sec}

 Manin~\cite{manin72} introduced the use of the continued fraction expansion of the rational $r$ to help speeding up the computation of $[r]^{\pm}$ considerably when the denominator of $r$ is large compared to $N$. See also~\cite{goldfeld} and Section~3.3.1 in~\cite{stein} for more details. However, we need to modify it slightly here as we should avoid non-unitary cusps if at all possible.

 \begin{defn}
 Recall that the set of right coset representatives of $\Gamma_0(N)$ in $\SL_2(\ZZ)$ is in bijection with $\PP^1\bigl(\cyclic{N}\bigr)$ by sending $\sm{a}{b}{c}{d}$ to $(c:d)$. For each such coset $\Gamma_0(N)\delta$ with $\delta= \sm{a}{b}{c}{d}$, we define the \textit{Manin symbol} by
  \begin{equation}\label{manin_symbol_eq}
    M(c:d) = 2\pi i \int_{i\infty}^0 f\vert_{\delta}(z) dz = 2\pi i \int_{a/c}^{b/d} f(z) dz= \lambda\Bigl(\frac{b}{d}\Bigr) - \lambda\Bigl(\frac{a}{c}\Bigr).
  \end{equation}
 \end{defn}
 We start by explaining how to reduce the computation of $[r]^{\pm}$ for a large denominator of $r$ to the computation of Manin symbols and then explain how to evaluate Manin symbols.

\subsection{Using continued fractions}\label{cont_frac_sec}

  Here is the original trick by Manin. We are given a rational number $r=a/m$. Consider the sequence of convergents of the continuous fraction of $r$:
  \begin{equation*}
    \frac{a_{-1}}{m_{-1}} = \frac{1}{0},\qquad
    \frac{a_0}{m_0} = \frac{a_0}{1},\qquad \dots,\qquad
    \frac{a_n}{m_n} =\frac{a}{m}.
  \end{equation*}
  We have $a_k m_{k-1} - a_{k-1}m_k = (-1)^{k-1}$. So the matrix
  \begin{equation*}
    \begin{pmatrix} a_k & (-1)^{k-1} a_{k-1} \\ m_k & (-1)^{k-1}m_{k-1}\end{pmatrix}
  \end{equation*}
  belongs to $\SL_2(\ZZ)$ and it sends any path linking $0$ to $i\infty$ to a path from $\frac{a_{k-1}}{m_{k-1}}$ to $\frac{a_k}{m_k}$. We find
  \begin{align*}
    \lambda(r) &= -2\pi i\cdot \Biggl( \int_{a_n/m_n}^{a_{n-1}/m_{n-1}} + \int_{a_{n-1}/m_{n-1}}^{a_{n-2}/m_{n-2}}
    +\cdots + \int_{a_1/m_1}^{a_0/m_0}+\int_{a_0}^{i\infty}\Biggr) f(z)dz  \\
    &= M\bigl(m_n:(-1)^{n-1}m_{n-1}\bigr) + M\bigl(m_{n-1}:(-1)^{n-2}m_{n-2}\bigr)+\cdots + M\bigl(m_1:1\bigr) + M\bigl(1:0\bigr).
  \end{align*}
  This allows to compute $\lambda(r)$ as a sum of Manin symbols $M(c:d)$, each of which is a modular symbol between two rational numbers of denominator $c$ and $d$ smaller than $N$.

  Now the problem with this way of splitting up is the following: Even if $r$ is a unitary cusp, it may be that some intermediate convergent $a_k/m_k$ is not unitary. Here is an adaptation, which may take a few steps more, but tries to avoid non-unitary cusps. In the end this is a great gain of speed.
  \begin{algo}{Algorithm: Try to split up the path into unitary Manin symbols}
    \item[Initialisation] Given $r=a/m$. If $m=1$, return $\lambda(0)$.
    \item[Find new cusp] Compute with the extended euclidean algorithm $x$ and $y$ such that $a\, y + x\,m =1$. Make sure that $-m/2 < y \leq m/2$.
    \item[Unitary?] If $-x/y$ is unitary, set $r'=-x/y$. Otherwise, set $r'=(x+\sign(y)\,a)/(y-\sign(y)\,m)$ if that is unitary. If both are non-unitary, set $r'=-x/y$.
    \item[Recursion] Call this function recursively with $r'$ and add the result to the Manin symbol $M(m:y)$.
  \end{algo}
  Here is an example of a case when both choices of cusps are non-unitary: For $N=36$ and $r=\tfrac{2}{5}$, neither $\tfrac{1}{2}$ nor $\tfrac{1}{3}$ is unitary. This can only happen when the squares of two distinct primes divide $N$.

  Note that if we have to go for the second choice for the cusp, then we still have $\vert y \vert < m$, but not $\vert y \vert < m/2$. So we are not certain any more if the algorithm takes only $O(\log(m))$ steps. In practice, the algorithm is quite effective in avoiding non-unitary cusps. We tested all elliptic curves of conductor at most $1000$ which are not semistable and whose conductor cannot be decreased by a quadratic twist. Among all $a/m$ with $m<N$, there were 77\% such that the best choice for $r'$ is unitary, for 22\% the second best choice is unitary and only in 1.4\% we have to pass to a non-unitary cusp $r'$.

\subsection{Unitary Manin symbols} \label{manin_symbols_subsec}

  As explained above, we now have to compute the Manin symbol $M(c:d)$ as defined in~\eqref{manin_symbol_eq}. We assume here first that both $c$ and $d$ are denominators of unitary cusps. In this case, we say that the Manin symbol $M(c:d)$ is unitary.  Note, that once we computed $M(c:d)$, we also know $M(-d:c) = - M(c:d)$. This is the formula~(2.2.6) in~\cite{cremona}. Further $\overline{M(c:d)}=M(c:-d)$. Also, there is a three term relation $M(c:d) +M(c+d:-c)+M(d,-c-d)=0$; which can be used to compute a further value if two of them are known.

  There are now at least three possible ways of evaluating the Manin symbol $M(c:d)$. Either by direct or indirect integration or by using transportation. Further note that $M(c:d)$ only depends on $(c:d)$ in $\PP^1\bigl(\cyclic{N}\bigr)$ and we may improve the performance by choosing good representatives $c$ and $d$.

  First, in the cases when both $c$ and $d$ are coprime to $N$, we could transport them as both cusps $\tfrac{a}{c}$ and $\tfrac{b}{d}$ are $\Gamma_0(N)$-equivalent to $0$. From the Section~\ref{transportable_subsec}, we see that the speed will be at best equal to $\tfrac{1}{N}$ and hence this method will usually lose out on the others below.

  Let $Q$ be the width of $\tfrac{a}{c}$ and $Q'$ be the width of $\tfrac{b}{d}$. Then the speed of using the direct integration from $\tfrac{a}{c}$ to $\tfrac{b}{d}$ is equal to
  \begin{equation*}
    \frac{1}{\sqrt{QQ'} \cdot \vert ad-bc\vert } =\frac{1}{\sqrt{QQ'}} \geq \frac{1}{N}
  \end{equation*}
  as seen in Section~\ref{cusp_to_cusp_subsec}. In the (most frequent) case when $c$ and $d$ are coprime to $N$, then the speed is indeed equal to $\tfrac{1}{N}$. Neglecting the contribution from $\ve$, this means that we expect a single sum over approximately $\tfrac{1}{2\pi} N \log(N)$ terms.

  By Section~\ref{ioo_to_cusp_subsec}, the indirect integration via $i\infty$ instead uses two sums with speed $\bigl(\vert c\vert \sqrt{Q}\bigr)^{-1}$ and $\bigl(\vert d \vert \sqrt{Q'}\bigr)^{-1}$ each. If we neglect again the contribution from $\varepsilon$, we expect in the case $\gcd(cd,N)=1$ to sum in total about
  \begin{equation*}
    \tfrac{\sqrt{N}}{2\pi} \Bigl( \vert c \vert\, \log \bigl(\vert c\vert \sqrt{N}\bigr) +  \vert d \vert\, \log\bigl(\vert d\vert \sqrt{N}\bigr) \Bigr ).
  \end{equation*}
  In particular, if we can find $c$ and $d$ representing the point on the projective line with $\vert c\vert$ and $\vert d\vert$ both smaller than $\tfrac{1}{2}\sqrt{N}$, then the indirect method is faster. This leads to the problem of finding small $c$ and $d$. The following lemma shows that we may just as well try to minimise $\vert c\vert + \vert d \vert$.

  \begin{lem}\label{tlogt_lem}
    Let $C=\sqrt{N}/(2\pi)$. Let $\gamma\colon \RR^2\to \RR_{\geq 0}$ be the continuous function such that $\gamma(x,y) = C \bigl( \vert x\vert \log(\vert x\vert \sqrt{N}) + \vert y \vert \log(\vert y\vert\sqrt{N})\bigr)$ for $xy\neq 0$. Let $L\subset\ZZ^2$ be a set not containing the origin. Let $(x_0,y_0)$ be a point of $L$ at which $\gamma$ is minimal and let $(x_1,y_1)$ be a point in $L$ at which $\vert (x,y)\vert = \vert x\vert+\vert y\vert$ is minimal. Then
    \begin{equation*}
      \frac{\gamma(x_1,y_1)}{\gamma(x_0,y_0)} = 1 + \bigO\Bigl(\frac{1}{\log(N)}\Bigr).
    \end{equation*}
  \end{lem}
  \begin{proof}
    Write $A= \vert(x_1,y_1)\vert$. Since $\gamma$ is increasing on rays leaving from the origin, we see that
    \begin{align*}
      \gamma(x_1,y_1)& \leq \max\Bigl\{ \gamma(x,y) \Bigm\vert \vert (x,y)\vert = A \Bigr\}\\
      \gamma(x_0,y_0)& \geq \min\Bigl\{ \gamma(x,y) \Bigm\vert \vert (x,y)\vert = A \Bigr\}
    \end{align*}
    It is not hard to show that the maximum above is $C \cdot \log(A\sqrt{N})$ and the minimum is $C \cdot \log(A\sqrt{N}/2)$. Hence we find
    \begin{equation*}
      \frac{\gamma(x_1,y_1)}{\gamma(x_0,y_0)} \leq \frac{1}{1-\frac{\log(2)}{\log(A\sqrt{N})}} = 1 + \bigO\Bigl(\frac{1}{\log(N)}\Bigr).\qedhere
    \end{equation*}
  \end{proof}

\subsection{Small coordinates of projective points}\label{smallp1_sec}

  Let $N$ be an integer and $(u:v)\in\PP^1\bigl(\cyclic{N}\bigr)$. In the above computation of Manin symbols, we came across the problem of finding the integers $c$ and $d$ such that $(u:v) = (c:d)$ and $\vert c \vert+ \vert d \vert$ is as small as possible. Write $\vert (c, d)\vert = \vert c\vert+\vert d \vert$ and $\Vert (c,d)\Vert = \sqrt{c^2+d^2}$.

  We are looking for the smallest non-zero vector in the lattice
  \begin{equation*}
    \Lambda_{(u:v)} = \Bigl\{ (c,d)\in\mathbb{Z}^2 \ \Bigl\vert \ c\cdot v \equiv d\cdot u \pmod{N} \Bigr\}
  \end{equation*}
  such that $\gcd(c,d)=1$.  Here is an algorithm based on Algorithm 1.3.14 in \cite{cohen}.
  \begin{algo}{Algorithm: Find good representatives for projective points}
    \item[Initialise] Set $\vec{x},\vec{y}$ to be a $\ZZ$-basis of $\Lambda_{(u:v)}$. If one of the coordinates $u$ or $v$ is invertible modulo $N$, say $v$, then we can do this as follows: Set $w$ to be the product of $u$ and the inverse of $v$ modulo $N$. Let $\vec{x} = (1,w)$ and $\vec{y} = (0,N)$. In the general case, we set $p=\gcd(u,N)$ and $q =\gcd(v,N)$; note that they must be coprime. Set $w$ to be the product of $\tfrac{u}{p}$ and the inverse of $\frac{v}{q}$ modulo $\tfrac{N}{pq}$. Then $\vec{x} = (\tfrac{N}{q},0)$ and $\vec{y}=(w\cdot p,q)$ is a basis.
    \item[Euclidean step] If the signs of $x_0$ and $x_1$ agree, then set $r$ to be the greatest integer smaller than $\frac{y_0+y_1}{x_0+x_1}$. Otherwise set $r$ to be the greatest integer smaller than $\frac{y_0-y_1}{x_0-x_1}$. Set $\vec{z}=\vec{y}-r\cdot\vec{x}$. If $\vert{\vec{z}-\vec{x}}\vert < \vert\vec{z}\vert$, then replace $\vec{z}$ by $\vec{z}-\vec{x}$.
    \item[Finished ?] If $\vert{\vec{z}}\vert < \vert{\vec{x}}\vert$, then set $\vec{y}$ to $\vec{x}$ and $\vec{x}$ to $\vec{z}$ and go back to the second step. Otherwise we can terminate the algorithm. If the coordinates of $\vec{x}$ are coprime, we return $\vec{x}$. If not, we run through small linear combinations of $\vec{x}$ and $\vec{z}$, starting with $\vec{z}$, until we hit one with coprime coordinates.
  \end{algo}

  The proof is very analogous to the one in~\cite{cohen}. As long as we do the second step, we know that $\vec{x}$ and $\vec{y}$ are a $\ZZ$-basis of the lattice $\Lambda_{(u:v)}$. The integer $r$ is chosen such that $\vert{\vec{z}}\vert$ is minimal. At the stage when we terminate, we are certain that $\vec{x}$ is the shortest non-zero vector of the lattice and $\vec{z}$ is the shortest, which is not a multiple of $\vec{x}$. The convex body theorem of Minkowski applied to the set of vectors of $\vert\cdot\vert$-norm at most $\sqrt{2N}$ guarantees that $\vert{\vec{x}}\vert\leq \sqrt{2N}$.

  Unfortunately, we cannot be certain that the algorithm will return the best of all choices. For instance with $N=30$ and $(u:v) = (11:1)$, we find that the shortest non-zero vector is $\vec{x} = (3,3)$ and the second minimum is $\vec{z} = (5,-5)$. None of them is allowed to represent $(11:1)$ in $\PP^1\bigl(\cyclic{30}\bigr)$. Even $\vec{x} + \vec{z} = (8,-2)$ and $\vec{z}-\vec{x} = (2,8)$ are not permitted. Only when we compute $2\vec{x} - \vec{z} = (11,1)$ and $\vec{x} + 2\vec{z} = (13,-7)$ will we find coprime coordinates. It is now not certain that the algorithm will find the shorter one first. Note that in this example $\vert {(11,1)} \vert = 12$ is much larger than $\sqrt{60}$.

  The following is a theoretical result about small coordinates for projective points that will be used later in Section~\ref{complexity_all_manin_sec}.

  \begin{lem}\label{smallp1_lem}
    There exists an absolute constant $C$ with the following property. Let $N$ be a square-free integer and let $P=(u:v)\in \PP^1\bigl(\cyclic{N}\bigr)$. Let $\vec{v}_P$ be the shortest non-zero vector in $\Lambda_P$ and let $\vec{w}_P$ be the shortest vector in $\Lambda_P$ which is not collinear to $\vec{v}_P$. Then there exists $\lambda\in \ZZ$ with $\vert \lambda \vert \leq C\cdot \log(N)^2$ such that the two coordinates of $\vec{w}_P+\lambda\vec{v}_P$ are coprime.

    In particular, there exists $(c,d)$ such that $(c:d)=(u:v)$ and
    \begin{equation*}
     \max\bigl( \vert c\vert , \vert d\vert \bigr)\leq  \frac{N}{\Vert \vec{v}_P\Vert} + C' \log(N)^2 \Vert\vec{v}_P\Vert
    \end{equation*}
    for some absolute constant $C'$.
  \end{lem}

  Note that it is vain to hope for a better bound, for instance independent of the size of $\Vert\vec{v}_P\Vert$. Suppose $N=2n$ is even. Then the size of the coordinates of the point $P = (1:n)$ cannot be decreased. For this example $\vec{v}_P=(2,0)$ is very small.

  \begin{proof}
    We will call the content, written $\co(x,y)$, of a point $(x,y)$ in $\ZZ^2$ the greatest common divisor of the two coordinates $x$ and $y$. Since $(u,v)\in\Lambda_{P}$ and $u$ and $v$ are coprime, there exists at least one point with content $1$ in $\Lambda_P$. It follows that the contents of two basis vectors of $\Lambda_P$ must be coprime integers. In particular $\co(\vec{v}_P)$ and $\co(\vec{w}_P)$ are coprime.

    Let $\vec{z}=(x,y)$ be a vector in $\Lambda_P$. Then there exists an integer $k$ such that $cx-dy=kN$. If $b = \gcd\bigl(k, \co(\vec{z})\bigr)$, then $\bigl(x/b,y/b\bigr)$ also belongs to $\Lambda_P$. Hence if we assume now that $\vec{z}$ is not divisible by any integer greater than $1$, then $b=1$. Thus $\co(\vec{z})$ divides $N$. In particular all points on the line $\mathcal{L} = \bigl\{ \vec{w}_P + \lambda \vec{v}_P \bigm\vert \lambda\in \ZZ\bigr\}$ have contents equal to a divisor of $N$.

    Consider two points $\vec{z}= \vec{w}_P + \lambda \vec{v}_P$ and $\vec{z}'=  \vec{w}_P + \lambda' \vec{v}_P$ on the line $\mathcal{L}$. We claim that the greatest common divisor of $\co(\vec{z})$ and $\co(\vec{z}')$ divides $\lambda-\lambda'$: It is not hard to show that this greatest common divisor divides $(\lambda'-\lambda)\cdot \gcd\bigl(\co(\vec{v}_P),\co(\vec{w}_P)\bigr)$ and so the above justifies the claim.

    For each prime divisor $\ell\mid N$, either $\ell$ does not divide the content of any point on $\mathcal{L}$ or the content of every $\ell$-th point is divisible by $\ell$. Let $\tilde N$ be the product of the prime divisors of $N$ dividing the content of one of the points on $\mathcal{L}$. The sequence $\co(\vec{w}_P+\lambda\vec{v}_P)$ as $\lambda$ varies in $\ZZ$ is periodic with period $\tilde N$. There is $\lambda_0$ such that $\vec{z}_0=\vec{w}_P+\lambda_0\vec{v}_P\in\mathcal{L}$ has content $\tilde N$. Now the content of $\vec{w}_P+\lambda\vec{v}_P$ is $\gcd(\lambda-\lambda_0,\tilde N)$.

    By a theorem of Iwaniec~\cite{iwaniec} on the Jacobsthal function, there is a constant $C$ such that any set of $C(\log(N))^2$ consecutive integers contain at least a unit modulo $N$. It follows that in the set $X\subset\mathcal{L}$ of $\vec{w}_P+\lambda \vec{v}_P$ with $\vert \lambda\vert \leq C/2 \log(N)^2 $ there is a point whose coordinates are coprime.

    \begin{figure}[hb]
   \centering
   \includegraphics[width=65mm,height=65mm]{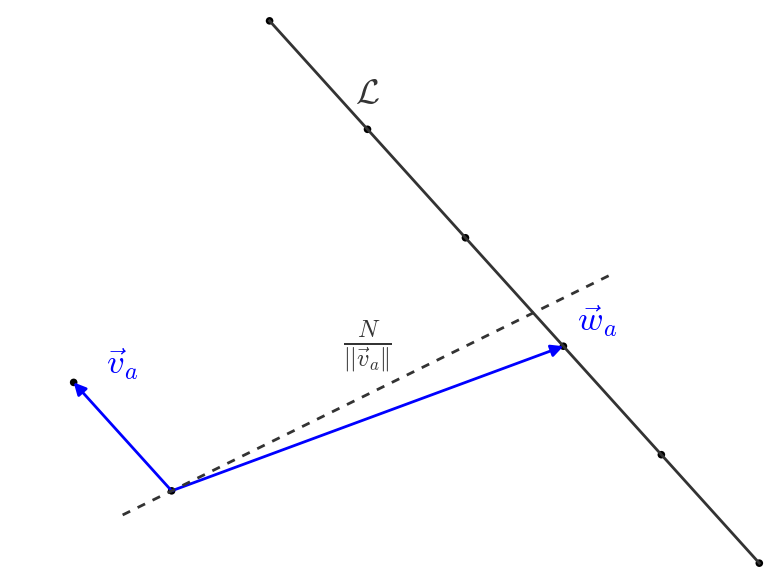}
   \caption{The two shortest vectors and the line $\mathcal{L}$}\label{p1_fig}
  \end{figure}

    The last sentence of the lemma follows from geometric considerations (see Figure~\ref{p1_fig}) measuring the length of this vector in $X$: The distance from $(0,0)$ to the real line containing $\mathcal{L}$ is $N/\Vert \vec{v}_P\Vert$. The length of the point $(x,y)$ in the set $X$ furthest away from $(0,0)$ satisfies
    \begin{equation*}
      \Vert (x,y)\Vert \leq \frac{N}{\Vert \vec{v}_P\Vert} + \bigl( C/2 \log(N)^2 + 1\bigr) \Vert \vec{v}_P\Vert
    \end{equation*}
    by the triangle inequality. Finally we use $\max(\vert x \vert , \vert y \vert ) \leq  \Vert(x,y)\Vert$.
  \end{proof}
   We also remark that when $N=p$ is prime, we have the much better bound $\vert c\vert + \vert d\vert \leq \sqrt{2N}$: The content of $\vec{v}_P$ can only be $1$ or $p$. But if it were $p$, then the representation of the form $P=(1:d)$ with $0\leq d <p$ or $(0:1)$ would be a smaller vector in $\Lambda_P$. Hence the shortest vector is always the best way to represent the point on $\PP^1\bigl(\cyclic{p}\bigr)$. By Minkowski's convex body theorem $\vert\vec{v}_P\vert \leq \sqrt{2N}$.

\subsection{Non-unitary Manin symbols}
  Let $(c:d)$ be such that at least one of them is not the denominator of a unitary cusp. For simplicity, we assume that $d$ is the denominator of a unitary cusp and $c$ is not. Given how much harder it is to work with non-unitary cusps, we should compute $M(c:d)$ as $\lambda(b/d)-\lambda(a/c)$ and we have to make $c$ as small as possible.

  Given an integer $N$ and $(u:v)\in\PP^1\bigl(\cyclic{N}\bigr)$, we are looking for $(c:d)=(u:v)$ such that $\vert c\vert$ is minimal. Let $M=(u,N)$ and $Q$ such that $N=MQ$. We can take $c=M$, which is minimal. The other coordinate $d$ has now to satisfy $Mv\equiv du\pmod{N}$ and $(M,d)=1$. Let $x$ and $y$ such that $xu + yN =M$. The congruence condition becomes $d\equiv xv\pmod{Q}$. Our first choice would be to take $d=xv$. However in case $xv$ and $M$ are not coprime, we add $Q$ to $xv$ until it  becomes coprime to $M$.


\section{Tweaks}\label{tweaks_sec}

  In this section, we present two ideas to make certain computations faster.

\subsection{Using partial sums} \label{partial_sec}

  Let $E$ be an elliptic curve over $\QQ$. Let $m$ be a small positive integer. Here is an idea that is useful for the evaluation of all symbols $[\tfrac{a}{m}]^{\pm}$ as $a$ varies through all integers $1\leq a < m$ coprime to $m$. In the application where we wish to evaluate a $p$-adic $L$-series for some small prime $p$, we would typically need this for $m=p^2$ or $p^3$. For the sake of simplicity we assume that $\tfrac{1}{m}$ is unitary.

  In equation~\eqref{sum_kappa_eq}, we have defined the partial sums
  \begin{equation*}
    \kappa_{j,m} (y) = \sum_{\substack{n\geq 1\\ n\equiv j \bmod{m}}} \frac{a_n}{n} \exp(-2\pi n y).
  \end{equation*}
  We have seen that we only need $m$ more terms in the sum to evaluate to a given precision all these partial sums for $j=0,\dots,m-1$.

  These can be used to evaluate $\lambda(\tau)$ whenever the real part $x$ of $\tau$ is a rational number with denominator $m$, say $x=\tfrac{a}{m}$:
  \begin{equation}\label{lambda_as_kappa_eq}
    \lambda\Bigl(\frac{a}{m} + yi\Bigr) = \sum_{j=0}^{m-1} \kappa_{j,m} (y) \cdot \zeta^{ja}
  \end{equation}
  where $\zeta = \exp(2\pi i/m)$. In case we are only interested in the plus modular symbols $[\cdot]^{+}$, we can do the computations with real numbers only.
  \begin{equation*}
    \re\Bigl( \lambda \bigl(\tfrac{a}{m} + yi\bigr)\Bigr) = \sum_{j=0}^{m-1} \kappa_{j,m} (y)\cdot \cos\bigl( 2 \pi a j / m \bigr).
  \end{equation*}
  We see here that it is possible to use fast Fourier transform if we are interested in evaluating $\lambda\bigl(\frac{a}{m}+yi\bigr)$ for all $a$ with a fixed $m$ and $y>0$. Note that the radix $m$ cannot be chosen to be a power of two, so we rely on mixed-radix algorithms. This has not yet been implemented in~\cite{tracticket}.

  We can use the above formula~\eqref{lambda_as_kappa_eq} together with equations~\eqref{sum_lambda_eq}, \eqref{split_r_ioo_eq} and~\eqref{Wrtau_eq}, to give a formula for the computation of $\lambda(\tfrac{a}{m})$ for all $a$ at once:
  \begin{equation}\label{lambda_rioo_eq}
    \lambda\Bigl(\frac{a}{m}\Bigr) = \sum_{n=1}^{\infty} \frac{a_n}{n} \cdot e^{-\frac{2\pi n}{m\sqrt{Q}}}\cdot \Bigl( e^{\frac{2\pi n a}{m} i} - \epsilon_Q e^{-\frac{2\pi n u}{m} i}\Bigr) = \sum_{j=0}^{m-1} \kappa_{j,m}\Bigl(\frac{1}{m\sqrt{Q}}\Bigr) \cdot \Bigl(\zeta_m^{ja} -\epsilon_Q \zeta_m^{-ju}\Bigr).
  \end{equation}
  where $u$ is an inverse of $Qa$ modulo $m$ and $\zeta_m = \exp(2\pi i /m)$.


  Similarly, we can express the direct integration from $r'=\tfrac{a'}{m'}$ to $r=\tfrac{a}{m}$ as a finite sum of partial sums: Let $Q$ and $Q'$ be the widths and set $d = \lcm(Q,Q')\cdot \vert am'-a'm\vert$ and $y = \sqrt{QQ'}\cdot \vert am'-a'm\vert$ and let $\tau$ be the optimal place in the upper half plan to cut the path in two, which we found in Section~\ref{cusp_to_cusp_subsec}.
  Then
  \begin{equation*}
    W_r(\tau)= \frac{\xi}{d} + \frac{i}{y}
  \end{equation*}
  where $\xi = (Qa'u+vm')Q'/\gcd(Q,Q')$ and $Qau +vm = 1$. Hence we obtain
  \begin{equation*}
    \lambda\bigl(\{r'\to r\}\bigr) = \sum_{j=0}^{d-1} \kappa_{j,d}\Bigl(\frac{1}{y}\Bigr) \Bigl( \epsilon_{Q'} \zeta_d^{j \xi'} - \epsilon_{Q} \zeta_d^{j \xi}\Bigr).
  \end{equation*}
  with $\xi'=(Q'au'+v'm)Q/\gcd(Q,Q')$ and $Q'a'u' +v'm' = 1$.
  Note however, that to use this formula only makes sense when $d$ is much smaller than $N$.

  Finally, we could also compute the transportable symbols using partial sums. For $\gamma=\bigl(\begin{smallmatrix} a& b\\ c &d\end{smallmatrix}\bigr)$, we find
  \begin{equation*}
    \lambda\bigl(\{r'\to r\}\bigr) = \sum_{j=0}^{\vert c\vert -1} \kappa_{j,\vert c\vert}\Bigl(\frac{1}{\vert c\vert}\Bigr) \Bigl(\zeta_c^{-dj} - \zeta_c^{aj}\Bigr).
  \end{equation*}

  We explain why it can be beneficial to use these partial sums: Even when computing a single one of these expressions, say $\lambda(\tfrac{a}{m})$ for some value of $m$, it may be worth wasting a bit of time and using the above formulae. We first compute all $\kappa_{j,m}(y)$ in one sum with $T(y,\ve)+m$ terms. Then we do one sum involving $m$ terms again. Hence if $m$ is small, say $m\ll \sqrt{N}$, we lose only very little time. Since $1/y^2$ is an integer in all cases above, it is easy to cache the values $\kappa_{j,m}(y)$ for later use. If we then encounter later another symbol with the same denominator $m$, we have to sum up only $m$ precomputed terms.

  However note that this is not practical for transportable symbols or for the computation of all Manin symbols as $m$ will be in the order of $N$ rather frequently.

\subsection{Quadratic twists}\label{twists_subsec}

  If $N$ is not square-free then one can often find a quadratic twist of the elliptic curve with smaller conductor. Since all the previous computations depend heavily on the conductor, it may be an advantage to do the computation on the twisted curve instead.

  Let $D$ be a fundamental discriminant such that the quadratic twist $E^{\dagger}$ by $D$ has minimal conductor among all quadratic twists of $E$. This needs not be unique, but for our considerations it does not seem to matter much which among them we choose. In practice we take the one with the largest period is best.

  Write $\sqrt{D}$ for the square root of $D$ in $\RR_{>0}$, if $D$ is positive, and in $i\RR_{>0}$, if $D$ is negative. We will use formula~(I.8.5) in~\cite{mtt}
  \begin{equation*}
    \lambda(r) = \frac{1}{\sqrt{D}} \sum_{u=1}^{\vert D \vert -1} \Bigl( \frac{D}{u} \Bigr) \lambda^{\dagger} \Bigl(r + \frac{u}{\vert D \vert}\Bigr)
  \end{equation*}
  where $\lambda^{\dagger}$ designates the modular symbol for the twisted elliptic curve.
  Since the rational numbers $r\pm \tfrac{u}{D}$ all have the same denominator, we can use the idea from the previous section to compute this sum with a single summation. Similar, if we wish to compute all modular symbols for $E$ with a given denominator.

  Note however that there is a small issue with this. Suppose $\ell$ is a prime dividing $D$ such that the conductor $N^{\dagger}$ of the twisted curve $E^{\dagger}$ is still divisible by $\ell^2$. This can happen for instance with $N=80$, $D=-4$, and $N^{\dagger}=40$. Now in this situation, we will evaluate modular symbols with denominator divisible by $\ell$. If $\ell$ did not divide the denominator of $r$, then the resulting cusp $r+\tfrac{u}{D}$ will not be unitary. Because our method is very much slower for non-unitary cusps, it is much better to avoid this. Hence we will remove all factors of $\ell$ in the fundamental discriminant if the twisted curve will still have additive reduction at $\ell$. Of course this affects only $\ell=2$ or $3$.

  How much do we expect this to speed up our computations?  We will use the notation $\bigO(f(N))$ to mean that the number of steps needed in the computations is, for sufficiently big $N$, bounded by $C\cdot f(N)$ for some constant $C>0$. Suppose we wish to evaluate $\lambda(r)$ for a rational $r$ with denominator $m$, which we suppose for simplicity to be coprime to $N$. We will compute about $\log(m)$ Manin symbols each with at worst a speed of $\tfrac{1}{N}$. So we will be summing about $\bigO\bigl(\log(m) \log(N) N\bigr)$ terms in total.

  Instead, using the twist by $D$, we will have $D$ times as many terms with a denominator of $m\cdot D$, but the conductor will be divided by $D\cdot D'$ where $D'$ is a factor of $D$. Hence we get about $\bigO\bigl( \log(mD) \log(N/DD') N/ D'\bigr)$ terms to sum. If $D'>1$, this is obviously a very good improvement. Otherwise it is negligible.

  The other major advantage of twisting is that there will be less non-unitary cusps on the twist. In particular when $E^{\dagger}$ is semistable, then  all cusps are unitary for $E^{\dagger}$. This way, we can compute even the non-unitary symbols for $E$ very quickly.


\section{Complexity}\label{complexity_sec}

  In~\cite{goldfeld} Goldfeld finds the complexity of evaluating one modular symbol on a semistable curve. We will refine this here. We will continue to use the notation $\bigO(f(N))$ introduced above to find an upper bound on the number of steps in an algorithm. Further the notation $\tilde \bigO(f(N))$ suppresses the possible further factors which are logarithmic in $f(N)$. As mentioned before, we neglect the issues with precision and simply find asymptotics for the number of terms that need to be summed up.

  We will assume throughout this section that $N$ is square-free; except for Section~\ref{complexity_all_den_subsec}. Recall that this implies that $E$ is semistable and hence the Manin constant $c_0$ is either $1$ or $2$ and hence the assumption made in Section~\ref{denominator_sec} can be neglected in this section.

\subsection{Periods}

  Although we have often neglected the size of $\ve$ in the previous consideration, we should find a proven lower bound for the size of the periods $\Omega^+$ and $\Omega^{-}$. This seems however difficult and the issue is already discussed in~\cite{goldfeld}.

\begin{figure}[ht]
   \centering
   \includegraphics[width=11cm]{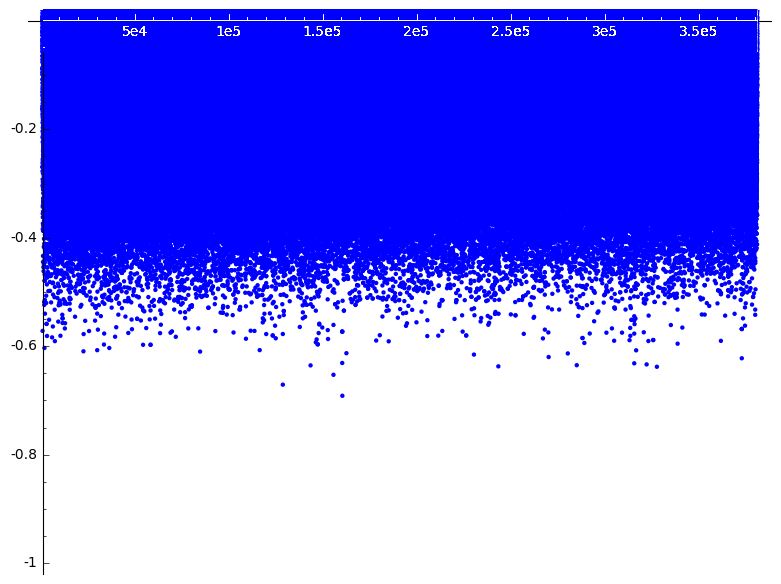}
   \caption{For each elliptic curve in the Cremona tables, the value of $\log\Omega^{+}/\log N$ on the vertical axis is compared with $N$ on the horizontal axis. Only negative values are plotted.}\label{period_fig}
  \end{figure}

  \begin{conj}[Goldfeld's period conjecture]
   There is a constant $\kappa>0$ such that $\Omega^+$ and $\Omega^-$ are larger than $\bigO\bigl(N^{-\kappa}\bigr)$ as $N\to\infty$.
  \end{conj}
  The graph in Figure~\ref{period_fig} presents numerical evidence in favour of this conjecture. In fact it looks like $\kappa=1$ is a very reasonable guess, while $\kappa<1/2$ is not likely.
  In Section~\ref{denominator_sec}, we have shown that for semistable curves the bound on the denominator of $[r]^{\pm}$ is at most $24$ for the strong Weil curve. Since the number of isogenous curves is also bounded, the denominator won't contribute to the asymptotic size of the error $\ve$. Under the conjecture above, we find that $-\log(\ve) = \bigO(\log(N))$.

  Without assuming the conjecture, it seems that one only knows (see~\cite{goldfeld}) that the periods are bounded by $\bigO\bigl(N^{-N}\bigr)$. This then gives a proven bound $-\log\ve =\bigO(N\log(N))$.

 \subsection{Fourier coefficients}

  We have to compute the coefficients $a_n$ for $n$ up to a bound $T$. In practice this is done by the command \texttt{ellan} in \texttt{PARI}. This function first computes the values $a_p$ for all primes up to $T$. When $p$ gets large, the preferred choice of algorithm for the Frobenius trace $a_p$ is the Schoof-Elkies-Atkin algorithm, which is known to run in polynomial time, with a heuristic expectation of $\tilde\bigO(\log^4p)$.  Hence to find all $a_p$ for $p<T$, we expect $\tilde \bigO(T)$ operations.
  The algorithm then uses the recursive formulae and the multiplicativity of $a_n$. This is done also in about $T$ steps. Therefore in total we expect $\tilde\bigO(T)$ operations.

  It is to be noted that in our implementation, this step does indeed take up a certain non-negligible portion of the total computation time. Initially, we precompute the first thousand coefficients $a_n$. If we later need more terms, we add them. However the way we interact with \texttt{PARI} currently it is faster to recompute all values from scratch unless we only have to add a small percentage of new values. Hence in practice, we may have to perform these computations more than once. For the theoretical considerations below, we may assume that we can determine  beforehand the highest value of $n$ ever needed and compute all values $a_n$ only once.

\subsection{Computing one modular symbol}

  Suppose $r=\frac{a}{m}\in \QQ$ with $0<a<m$ and we wish to evaluate $[r]^{\pm}$.
  As we supposed that $N$ is square-free, the cusp $r$ is unitary. We have seen in equations~\eqref{split_r_ioo_eq} and~\eqref{Wrtau_eq} that we can compute them by integrating to $\tau$ with imaginary part equal to $1/\bigl(m\sqrt{Q}\bigr)$ where $Q$ is the width of $r$. Lemma~\ref{truncate_lem} then gives us that we have to sum $T= \bigO\bigl(m\sqrt{Q}\log(m\sqrt{Q})\bigr) + \bigO(-\log(\ve)\,m \sqrt{Q})$ terms. For this we need to evaluate that many Fourier coefficients, but that is done in $\tilde\bigO(T)$ steps. As $Q\leq N$, we find that the total number of steps in the computation is $\tilde \bigO(m\sqrt{N})$ assuming Goldfeld's period conjecture.

  Of course, when $m$ is large, one should use Manin's trick in Section~\ref{manins_trick_sec} instead. Since $N$ is square-free, all cusps are unitary and hence we can split up the computation of $[r]^{\pm}$ into $\bigO(\log(m))$ Manin symbols. Now using the direct integration from cusp to cusp, any unitary Manin symbol can be computed in $\tilde\bigO(N\log(N))=\tilde\bigO(N)$ steps. We have now recovered
  \begin{thm}[Goldfeld, Theorem~2 in~\cite{goldfeld}]\label{goldfeld_thm}
    Assume Goldfeld's period conjecture holds. Then the modular symbol $\bigl[\tfrac{a}{m}\bigr]^{\pm}$ on a semistable curve $E$ defined over $\QQ$ of conductor $N$ can be computed in less than $\tilde\bigO(N\log(m))$ steps.
  \end{thm}
  However, we can often do much better. For instance, when $N$ is prime, then each Manin symbol can be computed in $\tilde \bigO(\sqrt{N}\log(N))=\tilde\bigO(\sqrt{N})$ steps due to the fact that projective coordinates can always be chosen of size $\bigO(\sqrt{N})$, see the remark after Lemma~\ref{smallp1_lem}. In fact a large proportion of Manin symbols are computable at that complexity:
  \begin{prop}
     Assume Goldfeld's period conjecture holds. For each $N$, there is a subset $\mathcal{P}$ containing at least $95\%$ of all points on $\PP^1\bigl(\ZZ/N\ZZ\bigr)$ such that each  Manin symbol $M(x)$ for $x\in\mathcal{P}$ can be computed in less than $\tilde\bigO(N^{1/2})$ steps.
  \end{prop}
  \begin{proof}
    Let $\vec v$ be a vector with $\Vert \vec{v}\Vert<\sqrt{N}$ and whose coordinates are coprime. Then $\vec{v}$ is the shortest vector in a lattice $\Lambda_P$ for some $P\in \PP^1\bigl(\cyclic{N}\bigr)$. Take $\mathcal{P}$ to be the set of all these points. For each $P$, there is only one other non-zero element of $\Lambda_P$, namely $-\vec{v}$, in the ball of radius $\sqrt{N}$. There are approximatively $\frac{1}{2} \frac{6}{\pi^2} \pi \sqrt{N}^2=\frac{3}{\pi}N$ pairs of opposite points with coprime integers in this ball. This is asymptotically more than $95\%$ of all elements in $\PP^1\bigl(\cyclic{N}\bigr)$.
  \end{proof}

 \subsection{Computing all modular symbol with a given small denominator}\label{complexity_all_den_subsec}

  Recall that $\delta^2$ is the largest square dividing $N$.
  \begin{thm}\label{complexity_allden_thm}
   Assume Goldfeld's period conjecture and assume Manin's conjecture that $c_0=1$. Let $m>1$ be an integer. Then there is a method to evaluate all modular symbols
   \begin{equation*}
     \Bigl\{\bigl[\tfrac{a}{m}\bigr]^{\pm}\ \Bigm\vert\  0\leq a <m \text{ and }\gcd(a,m)=1 \Bigr\}
   \end{equation*}
   for any elliptic curve over $\QQ$ of conductor $N$ with $\gcd(m,\delta)=1$ in less than $\tilde\bigO(N^{1/2})$ steps.
  \end{thm}
  If we restrict to semistable curves, the condition on $c_0$ can be dropped and $m$ is always coprime to $\delta=1$.
  \begin{proof}
   By assumption all cusps $a/m$ are unitary. Recall from the explanations in Section~\ref{partial_sec} the formula~\eqref{lambda_rioo_eq}. Hence we start by evaluating all $\bigl\{\kappa_{j,m}(y)\bigr\}_j$ with $y=m\sqrt{Q}$ using the approximation in Lemma~\ref{Tprime_lem}. This can be done with $m$ sums of $\bigO(m\sqrt{Q}/m)$ terms. Thus this first part takes $\bigO(m\sqrt{Q})$ steps.

   Given the vector $\bigl\{\kappa_{j,m}(1/y)\bigr\}_j$, we need to obtain the vector
   \begin{equation*}
     \Bigl\{\sum_{j=0}^{m-1} \kappa_{j,m}\bigl(\frac{1}{m\sqrt{Q}}\bigr) \zeta^{ja} \Bigm\vert a=1,\dots, m-1\Bigr\}
   \end{equation*}
   where $\zeta= \exp(2\pi i/m)$. For this we can use fast Fourier transform; in particular with Bluestein's multi-radix algorithm~\cite{bluestein} this is done in $\bigO(m\log(m))$ steps even when $m$ is not a prime power. Hence we get a complexity of $\tilde\bigO(m\sqrt{N})$ as $Q\leq N$, which yields the result as $m$ is fixed.
  \end{proof}

  In practice, we may be interested in computing approximations to the $p$-adic $L$-function for varying elliptic curves. Let $p^r$ be a fixed prime power. In order to determine the $r$-th approximation to the $p$-adic $L$-function as explained in~\cite{steinwuthrich}, we will only need to compute all modular symbols with denominator $p^r$. By the above this can be done with a complexity $\tilde\bigO(\sqrt{N})$.
  \begin{figure}[ht]
   \centering
   \includegraphics[width=11cm]{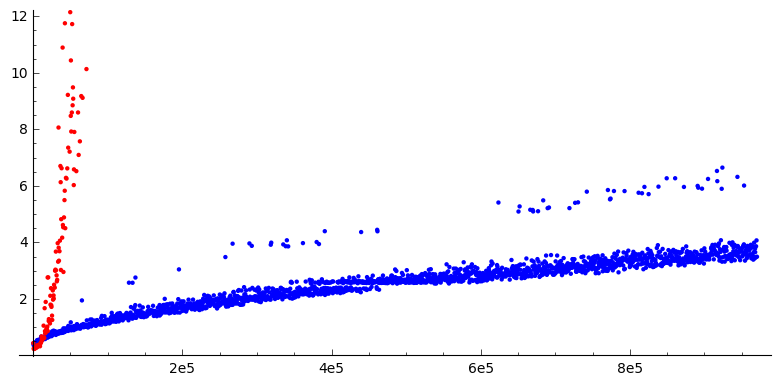}
   \caption{Comparison of approximations computing $5$-adic $L$-functions for semistable elliptic curves}\label{5adicl_fig}
  \end{figure}

  In Figure~\ref{5adicl_fig}, we plot the time to compute the fourth approximation $P_4$ in the notation of~\cite{steinwuthrich}. We tested random semistable curves with good ordinary reduction at~$5$ of conductor up to $10^6$, either from Cremona's table or from table of Stein and Watkins. The steeply increasing set of values uses \texttt{eclib}, the other timings are obtained with our implementation. The graph shows two anomalies: First there are a small number of values significantly higher than others. It turns out these are those examples for which the standard double precision of $53$~bits is not sufficient and the implementation has to use the much slower library of arbitrary precision floating point numbers. Secondly, there is a strange vertical strip empty. This is due to the choices of the values of $B(\varsigma)$ in~\eqref{Bvarsigma_eq}; these particular computations involve about $277200$~terms in the sum.

 \subsection{Computing all Manin symbols}\label{complexity_all_manin_sec}

  We wish to compare the numerical modular symbols to current implementations. Traditional methods start by finding a basis for the space of modular symbols attached to $E$ in the space of all modular symbols for $\Gamma_0(N)$. This is equivalent to computing all Manin symbols $M(c:d)$ for $(c:d)\in\PP^1\bigl(\cyclic{N}\bigr)$. We will estimate therefore the complexity to compute all Manin symbols via numerical approximations. Note however that in practice, we never do this. Instead we fill up the cached values for Manin symbols as we go along.
  \begin{thm}\label{complexity_thm}
   Assume Goldfeld's period conjecture is true. Then there is a method to evaluate all Manin symbols for any semistable elliptic curve over $\QQ$ of conductor $N$ in less than $\tilde\bigO(N^{7/4})$ steps.
  \end{thm}
  \begin{proof}
    As in Section~\ref{smallp1_sec} we denote for each $P\in \PP^1\bigl(\cyclic{N}\bigr)$ the lattice $\Lambda_P$ whose points with coprime coordinates are the possible representations of $P$. Let $\vec{v}_P$ be the shortest non-zero vector in $\Lambda_P$.

    We start by evaluating all $M(P)$ for those $P\in\PP^1\bigl(\cyclic{N}\bigr)$ with $\no{\vec{v}_P}\leq 2 N^{3/8}$. There are at most $4\pi N^{3/4}$ of them and each such Manin symbols can be evaluated in $\tilde\bigO(N)$ steps using the direct method. Hence all of them are done in $\tilde\bigO(N^{7/4})$ steps.

    Now, we may assume that $\no{\vec{v}_P}> 2 N^{3/8}$. By Minkowski's convex body theorem, we also know that $\no{\vec{v}_P}\leq 2/\sqrt{\pi}\, N^{1/2}$. We apply Lemma~\ref{smallp1_lem} and find that $P$ can be written as $(c:d)$ with
    \begin{equation*}
     \max(\vert c\vert,\vert d\vert )\leq  \frac{N}{\no{\vec{v}_P}} + C'\log(N)^2\,\no{\vec{v}_P} < \tfrac{1}{2}\, N^{5/8} + C'\log(N)^2 2/\sqrt{\pi}\, N^{1/2} =\bigO(  N^{5/8} )
    \end{equation*}
    if $N$ is sufficiently large. Therefore, we can evaluate all the remaining Manin symbols using the indirect method if we can compute all $\lambda(\tfrac{a}{m})$ with $m< N^{5/8}$ and $0<a<m$. To do this, we use the idea in the previous section and we can get all $\lambda(\tfrac{a}{m})$ for a fixed $m$ in $\tilde\bigO(m\sqrt{N})$ steps. Hence to find all $\lambda(\tfrac{a}{m})$ for $m<N^{5/8}$ we require $\tilde\bigO\bigl((N^{5/8})^2\sqrt{N}\bigr) = \tilde\bigO(N^{7/4})$ steps.
  \end{proof}

  Again, we can comment that this complexity is not always optimal. If $N$ is prime, then all Manin symbols can be computed using all $\lambda(\tfrac{a}{m})$ with $m<\bigO(\sqrt{N})$. This gives a total complexity of $\tilde\bigO(N^{3/2})$.

  To get unconditional results, i.e., independent of Goldfeld's conjecture, one may multiply all the complexities above with $N$.

   The current implementations involve Gaussian elimination on sparse matrices of size $\bigO(N)\times\bigO(N)$. More precisely, as explained in Algorithm~8.38 in~\cite{stein}, each matrix has about $N/3$ rows each containing at most three non-zero values. It is not hard to see that Gaussian elimination needs at least $\bigO(N^{3/2})$ steps for each such matrix as we expect to reach a dense matrix by the time we are dealing with the last $\sqrt{N}$ rows. However it would be rather hard to prove a precise complexity for the full algorithm.


\section{Examples}\label{example_sec}

 The computations below are performed with our implementation~\cite{tracticket} written in \texttt{Cython}~\cite{cython}. Note that this implementation is not fully optimised. The emphasis was on getting correct results for unitary cusps and for computing all modular symbols for a given denominator. For instance, it does not include the algorithm with the complexity of Theorem~\ref{complexity_thm}, though for the range of considered conductors this will not matter much.

 First, we present a concrete example of our methods. We choose the curve $E=$ 234446a1, famous for being the first curve in Cremona's tables of rank $4$. It is semistable so we do not have to worry about non-unitary cusps. We are interested in computing the $p$-adic $L$-function $\mathscr{L}_p(E,T)$ as explained in~\cite{steinwuthrich} at the good ordinary prime $p=5$. There are no isogenies from $E$ defined over $\QQ$ and the N\'eron period lattice $\Lambda_E = 1.486336\dots \ZZ \oplus 0.800625\dots \ZZ i$ is rectangular. Therefore the modular symbols $[r]^{\pm}$ are integers. In fact $[\tfrac{1}{27}]^+=[\tfrac{1}{7}]^{-} = 1$ and $[\tfrac{1}{7}]^+=0$ show that the values $\lambda(r)$ generate $\Lambda_E$. In particular, we have to approximate the real part of $\lambda(r)$ to precision $0.743168$. When computing all values of $[\tfrac{a}{5}]^+$ using the partial sums $\kappa_{j,5}(y)$ we need $T=2923$ terms and the precision of $53$ bits is enough. The largest error in evaluating these was smaller than $0.00032\,\Omega^+$. Similar for all values $[\tfrac{a}{25}]^{+}$ we only need to sum $17716$ terms, still with precision of $53$ bits. Using these values one finds that the fourth coefficient of $\mathscr{L}_5(E,T)$ is congruent to $1$ modulo $5$. This implies that the rank of $E(\QQ)$ is at most $4$. Together with the explicit basis of $E(\QQ)$ one can deduce without much further effort that the $5$-primary part of the Tate-Shafarevich group $\Sha(E/\QQ)$ is trivial.

 Next, in comparison an example involving non-unitary cusps. Let $E$ be the elliptic curve 1017a1, which has additive reduction at $3$ of Kodaira type III. Its quadratic twist by $-3$ is 1017e1, which has type III${}^{*}$ at $3$. The seemingly harmless computation of $[\tfrac{1}{3}]^+$ now involves more than 48000 terms to sum in total. Instead $[\tfrac{1}{5}]^+$ only requires 217 terms to sum. Though we have to admit that it is likely that the implementation for the non-unitary cases could be improved.

\begin{figure}[h]
 \centering
 \includegraphics[width=12cm]{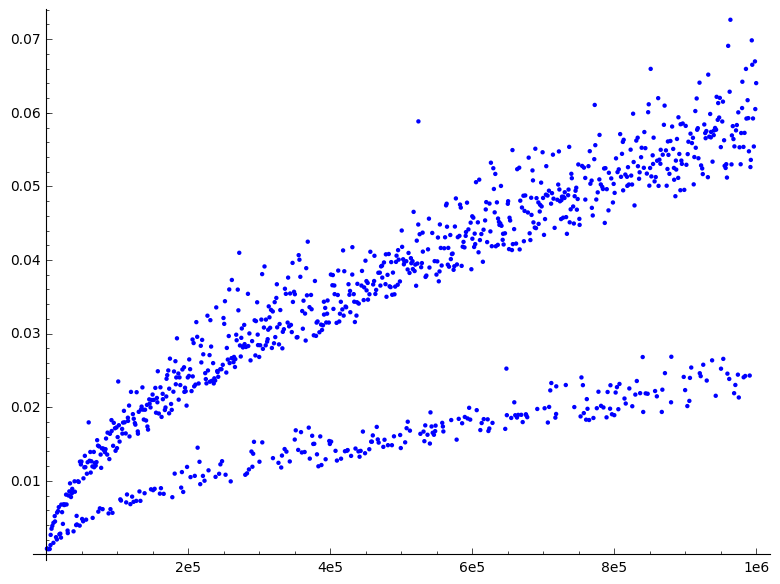}
 \caption{Time to compute all symbols $[\frac{a}{25}]^{+}$ for some semistable curves. }\label{den25_fig}
\end{figure}

 Now to the asymptotic behaviour as $N$ increases. In Figure~\ref{den25_fig}, we used the numerical implementation to compute all $[\tfrac{a}{25}]^+$ for various random semistable curves. The time in seconds is plotted against the conductor $N$. The quicker ones are those with conductor divisible by $5$.
 This and the following computations were performed on rather standard hardware, for instance on a Intel Xeon E5-2660 2.6 GHz virtual machine with two cores.

 We now pass to compare the various implementations. There is our implementation~\cite{tracticket} of numerical modular symbols written in \texttt{Cython}~\cite{cython} incorporated into \texttt{SageMath}~\cite{sage}, the implementation of \texttt{eclib}~\cite{eclib}, written in \texttt{C}, also accessible within \texttt{SageMath}, the pure \texttt{Python} implementation in \texttt{SageMath}, the implementation in \texttt{Magma}~\cite{magmamanual} and the implementation in \texttt{PARI}~\cite{pari}. First we will exclude the pure \texttt{Python} implementation in \texttt{SageMath} and the one in \texttt{PARI}, which is still under development, as they are both significantly slower then the other three. The fact that these four implementations of the same algorithm have such different timings explains why we cannot compare them directly: they are written in different languages. Also, we call them from within \texttt{SageMath} and the time \texttt{SageMath} spends to call the underlying code varies much. Instead we want to illustrate the asymptotic behaviour of the computation.

\begin{figure}[h]
 \centering
 \includegraphics[width=15cm]{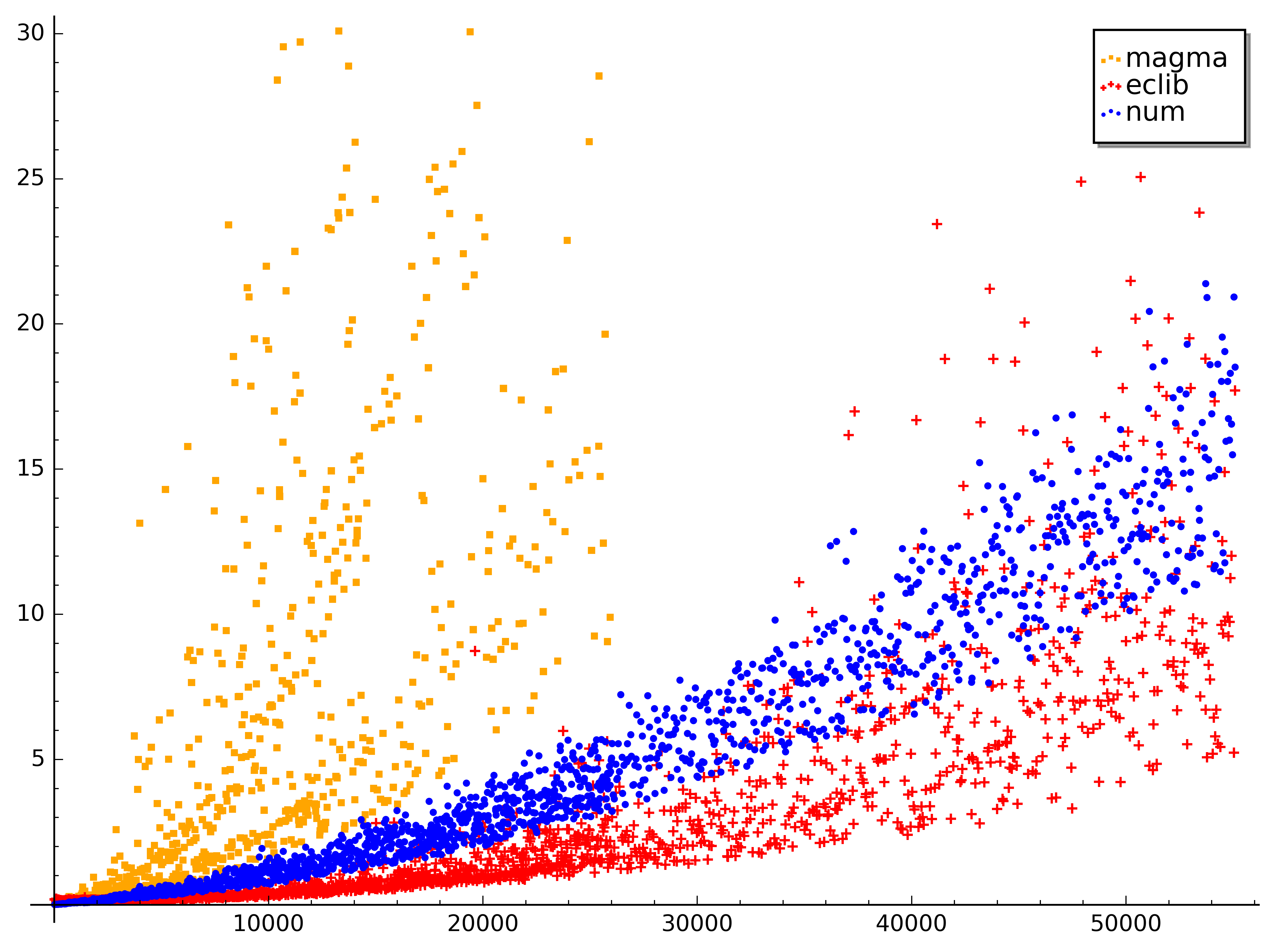}
 \caption{Time to compute all Manin symbols for some semistable curves}\label{all_fig}
\end{figure}

 In Figure~\ref{all_fig} we plot the time to compute all Manin symbols $M(c:d)$ using the numerical implementation (\textcolor{blue}{$\bullet$}) against the determination of the space of modular symbols by \texttt{Magma} (\textcolor{orange}{$\blacksquare$}) and \texttt{eclib} (\textcolor{red}{$\bm{+}$}). We do this in all three cases for random semistable curves of conductor up to $55000$. The computation was stopped after 30 seconds, meaning that for some curves the plotted point would lie an unknown amount above the visible part. The computations in \texttt{Magma} became rather quickly too complicated and they were stopped after conductor $25937$.


\bibliographystyle{amsplain}
\bibliography{modsym}

\end{document}